\RequirePackage{amsthm}
\documentclass[iicol,sn-mathphys-ay]{sn-jnl} 

\usepackage{graphicx}%
\usepackage{multirow}%
\usepackage{amsmath,amssymb,amsfonts}%
\usepackage{amsthm}%
\usepackage{mathrsfs}%
\usepackage[title]{appendix}%
\usepackage{xcolor}%
\usepackage{textcomp}%
\usepackage{manyfoot}%
\usepackage{booktabs}%
\usepackage{algorithm}%
\usepackage{algorithmicx}%
\usepackage{algpseudocode}%
\usepackage{listings}%
\usepackage{epstopdf}
\usepackage{caption}

\theoremstyle{thmstyleone}%
\newtheorem{theorem}{Theorem}
\theoremstyle{thmstyletwo}%
\newtheorem{example}{Example}%
\newtheorem{remark}{Remark}%

\theoremstyle{thmstylethree}%

\newtheorem{lemma}{Lemma}
\newtheorem{assumption}{Assumption}

\begin{document}

\title[Article Title]{Subsampling for Big Data Linear Models with Measurement Errors}


\author[1]{\fnm{Jiangshan} \sur{Ju}}\email{jsju1998@126.com}

\author[2]{\fnm{Min-Qian} \sur{Liu}}\email{mqliu@nankai.edu.cn}

\author*[1]{\fnm{Mingqiu} \sur{Wang}}\email{mqwang@vip.126.com}

\author[1]{\fnm{Shengli} \sur{Zhao}}\email{zhaoshli758@126.com}

\affil*[1]{\orgdiv{School of Statistics and Data Science}, \orgname{Qufu Normal University}, \orgaddress{\city{Qufu}, \postcode{273165},\country{China} }}

\affil[2]{\orgdiv{NITFID, LPMC \& KLMDASR, School of Statistics and Data Science}, \orgname{Nankai University}, \orgaddress{\city{Tianjin}, \postcode{300071},\country{China}}}


\abstract{Subsampling algorithms for various parametric regression models with massive data have been extensively investigated in recent years.
	However, all existing studies on subsampling heavily rely on clean massive data.
	In practical applications, the observed covariates may suffer from inaccuracies due to measurement errors.
	To address the challenge of large datasets with measurement errors, this study explores two subsampling algorithms based on the corrected likelihood approach: the optimal subsampling algorithm utilizing inverse probability weighting and the perturbation subsampling algorithm employing random weighting with a perfectly known distribution.
	Theoretical properties for both algorithms are provided.
	Numerical simulations and two real-data examples demonstrate the effectiveness of these proposed methods compared to other existing algorithms.}




\maketitle

\section{Introduction}\label{sec1}

To address the ever-increasing volume of data brought about by technological advancements, it is imperative to adopt refined techniques such as divide and conquer, online updating of streaming data, and subsampling-based methods.
These techniques offer effective solutions to computational challenges posed by large datasets.
However, existing literature often assumes direct and accurate observation of covariates, which may not always be feasible in practical data collection scenarios.
Consequently, statistical models that do not account for measurement errors can lead to biased estimation results.
Therefore, it is essential to investigate subsampling algorithms for linear models that take measurement errors in covariates into account.

For subsampling algorithms of linear models, \cite{MaMahoneyYu2015} proposed a leverage sampling algorithm based on leverage scores and their linear transformations.
\cite{Maetal2022} investigated the asymptotic normality and asymptotic unbiased properties of the subsample based estimator.
A deterministic subsampling method named information-based optimal subdata selection (IBOSS) was proposed by \cite{WangYang2019}, and extended by \cite{Wang2019b}, aiming to find subsamples with the maximum information matrix under the D-optimality criterion, which performs well in finding corners.
\cite{ChengWangYang2020} and \cite{YuLiuWang2023} extended the IBOSS algorithm to logistic and nonlinear models, respectively.
\cite{WangElmstedtWongXu2021} proposed an orthogonal subsampling approach for big data linear regression.
\cite{Chasiotis2025} evaluated the algorithm of the LEVSS approach \citep{Yu2022} for the selection of the most informative data points in order to estimate
unknown parameters.
\cite{MengZhang2020, MengYu2022}, \cite{YiZhou2023} and \cite{ZhangZhou2024} explored the method of using space-filling or uniform designs to obtain the subsample so that a wide range of models could be considered.

\cite{WangZhuMa2018} proposed the optimal subsampling method based on the A-optimality criterion. This method was further developed by \cite{Wang2019a}, \cite{AiWangYuZhang2021}, \cite{AiYuZhangWang2021}, \cite{WangMa2021}, and \cite{YuWangAiZhang2022}.
\cite{WangKim2022} introduced the maximum sample conditional likelihood estimation, and enhanced the estimator for selected subsamples. This approach overcomes the limitations of inverse probability weighting and makes more efficient use of sample information.
Due to the fact that inverse probability weighting requires calculation of the probabilities of all samples simultaneously, which imposes big computational burden, \cite{YaoJin2024} introduced a perturbation subsampling method. This method utilizes repeated random weighting based on known distributions to address the limitations of inverse probability weighting. It has been successfully applied to linear models, longitudinal data, and high-dimensional data, demonstrating promising performance.

For the measurement error model, \cite{Fuller1987} systematically introduced a comprehensive statistical inference of linear regression models with measurement errors.
\cite{Nakamura1990} proposed the corrected score method for the generalized linear model with measurement errors.
\cite{Liang1999} offered the parameter estimation of a semi-parametric partially linear model with measurement errors.
\cite{Carroll2006} systematically studied the theory of nonlinear regression models with measurement errors.
\cite{LiangLi2009} examined variable selection in partially linear models with measurement errors, and suggested that when the variance of the measurement error is unknown, it can be estimated through repeated observations.
\cite{LeeWangSchifano2020} introduced an online update method for correcting measurement errors in big data streams.

This study focuses on the subsampling problem in linear models with measurement errors. The presence of measurement errors in covariates can lead to inaccuracies of parameter estimation, thereby diminishing the statistical power. We employ the corrected likelihood approach proposed by \cite{Nakamura1990} to estimate the parameters with subsamples. An optimal subsampling method based on the corrected likelihood approach is proposed, and the optimal subsampling probabilities are determined by minimizing the trace of the variance. The consistency and asymptotic normality of estimators obtained by this approach are established. {\color{red}Although the theoretical results and numerical studies illustrate that our method perform well for big data with measurement errors, the computing times of the optimal subsampling method is the drawback. The reason is that the explicit solution does not need the iterative algorithm for linear models even when there are measurement errors in covariates.} Furthermore, we propose a perturbation subsampling method based on the corrected likelihood approach that approximates the objective function of the full data using a perturbation with independently generated stochastic weights. The effectiveness of the proposed methods is also confirmed through numerical analysis. By accounting for measurement errors in covariates, more precise and reliable results can be obtained in the analysis of massive datasets.
Our approaches not only alleviate the computational burden associated with parameter estimation in big data but also enhance computational efficiency and improve prediction accuracy.

The rest of the paper is outlined as follows.
Section 2 offers a comprehensive introduction to the model setup and parameter estimation of the measurement error model.
Sections 3 and 4 propose the linear model subsampling algorithms with measurement errors and establish the correspondingly theoretical properties.
Section 5 comprises numerical simulations. Section 6 presents case studies aimed at validating effectiveness of the algorithms.
Finally, Section 7 summarizes this paper. The detailed proofs are provided in Appendix A.

\section{Linear Model with Measurement Errors}\label{sec2}

In this section, we present an overview of the model and parameter estimation.

\subsection{Model}\label{subsec2}

Here, we consider the linear model with measurement errors
\begin{equation}\label{eq1}
	\left\{
	\begin{aligned}
		y_{i} = \mathbf{X}_{i}^{T} \boldsymbol{\beta} + \epsilon_{i}\\
		\mathbf{W}_{i} = \mathbf{X}_{i} + \mathbf{U}_{i}
	\end{aligned}
	\right. ,  ~i = 1 , 2 , \ldots , n,
\end{equation}
where $\mathbf{X}_{i} = (x_{i1},\ldots,x_{ip})^{T}$ is a $p$-dimensional covariate vector, $y_{i}$ is the corresponding response variable,
$\boldsymbol{\beta} = (\beta_{1}, \beta_{2},\ldots,\beta_{p})^{T}$ is an unknown parameter vector,
and $\epsilon_{i} $ is a random error term with mean zero and variance $\sigma^{2}$.
However, in applications, the exact value of $\mathbf{X}_{i}$ is often difficult to obtain.
Let $\mathbf{U}_{i} = (u_{i1},\ldots,u_{ip})^{T}$ be a random error vector with mean zero and variance-covariance matrix $\boldsymbol{\Sigma}_{uu}$,
and $\mathbf{W}_{i}$ be the actual observed random vector of $\mathbf{X}_{i}$.
Assuming $\boldsymbol{\Sigma}_{uu}$ is known, the measurement error $\mathbf{U}_{i}$ is independent of $\epsilon_{i}$ and $\mathbf{X}_{i}$. Denote the full data as $\mathcal{F}_{n} = \{(\mathbf{W}_{i}, y_{i}): i=1,\ldots,n\}$.

According to \cite{Fuller1987}, the ordinary least squares method cannot be directly applied to estimate linear models with measurement errors, as the resulting estimators are biased and inconsistent.
For model (\ref{eq1}), it follows that
\begin{equation*}
	y_{i} = \mathbf{W}_{i}^{T} \boldsymbol{\beta} + \epsilon_{i} - \mathbf{U}_{i}^{T} \boldsymbol{\beta} \triangleq \mathbf{W}_{i}^{T} \boldsymbol{\beta} + \delta_{i} ,~ i = 1 , 2 , \ldots , n ,
\end{equation*}
where $\delta_{i} = \epsilon_{i} - \mathbf{U}_{i}^{T} \boldsymbol{\beta}$. Note that
\begin{align*}
	Cov(\mathbf{W}_{i}, \delta_{i})
	&= Cov(\mathbf{X}_{i} + \mathbf{U}_{i}, \epsilon_{i} - \mathbf{U}_{i}^{T} \boldsymbol{\beta}) \\
	&= - \boldsymbol{\Sigma}_{uu} \boldsymbol{\beta}
	 \neq 0.
\end{align*}
Since the assumption of independence is violated, the ordinary method cannot be applied directly. This also suggests that when there are no measurement errors, i.e., $\boldsymbol{\Sigma}_{uu} = 0$, the estimator is unbiased. Example \ref{ex1} illustrates the influences of measurement errors.

\begin{figure*}[htb]
	\centering
	\begin{minipage}{5cm}
		\centering
		\includegraphics[width=5cm,height=5cm]{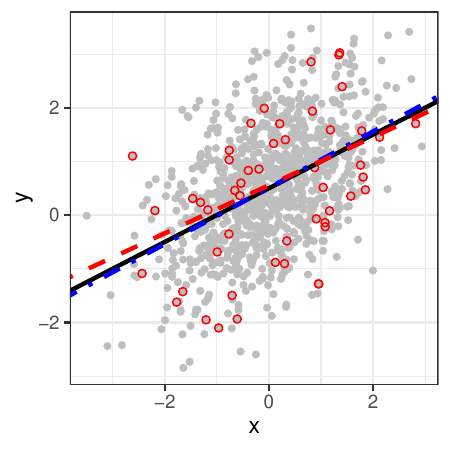}
		\caption*{(1)}
	\end{minipage}	
	\begin{minipage}{6cm}
		\centering
		\includegraphics[width=7.5cm,height=5cm]{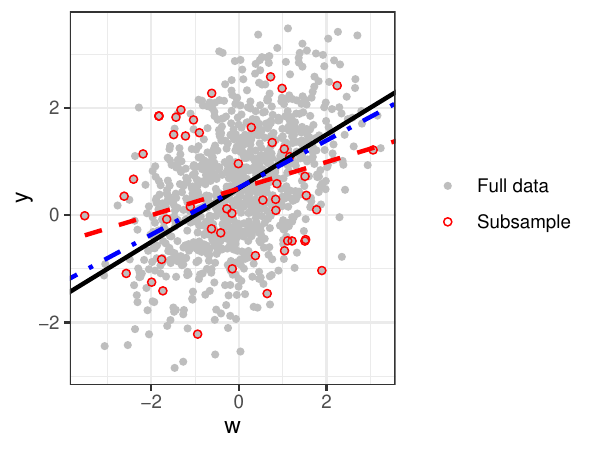}
		\caption*{(2)}
	\end{minipage}
	\caption{(1) The fitting results of the data without measurement errors. (2) The fitting results of the data with measurement errors, while ignoring measurement errors.}
	\label{fig1}
\end{figure*}

\begin{example}\label{ex1}
	Suppose the responses $y_i$ are generated from $y_i=0.5+0.5x_i+\varepsilon_i$,
	$x_i\sim N(0,1)$ and $\varepsilon_i\sim N(0,1)$, ~$i=1,\ldots,1000$.
	We consider the measurement error of $x_{i}$ by replacing $x_{i}$ with $w_i=x_i+u_i$, $u_i \sim N(0,0.5^{2})$.
	Figure \ref{fig1} describes the fitting results of the data with or without measurement errors.
	In Figure \ref{fig1}, the black solid line represents true model.
	The blue dotted line represents the regression line fitted based on the full dataset,
	while the red dashed line represents the regression line fitted {based on} a subsample with size 50 selected using L-optimal subsampling.
	From Figure \ref{fig1}, we can see that the presence of measurement errors in covariates have great effect on the resulting subsample estimator,
	if measurement errors are ignored.
\end{example}

\subsection{Parameter estimation}\label{subsec3}
To correct measurement errors, we apply the corrected likelihood method proposed by \cite{Nakamura1990}.
When $\boldsymbol{\Sigma}_{uu}$ is known, using the corrected likelihood method, we obtain
$$
\ell(\boldsymbol{\beta}) = \frac{1}{2n}\sum_{i=1}^{n}(y_{i}-\mathbf{W}_{i}^{T}\boldsymbol{\beta})^{2} - \frac{1}{2} \boldsymbol{\beta}^{T} \boldsymbol{\Sigma}_{uu} \boldsymbol{\beta}.
$$
Minimizing $\ell(\boldsymbol{\beta})$, we have
\begin{align*}
	\hat{\boldsymbol{\beta}} &= \arg\underset{\boldsymbol{\beta}}{\min} \ell(\boldsymbol{\beta}) \\
	&= \left(\sum_{i=1}^{n} \mathbf{W}_{i}\mathbf{W}_{i}^{T} -n\boldsymbol{\Sigma}_{uu} \right)^{-1} \sum_{i=1}^{n}\mathbf{W}_{i}y_{i}.
\end{align*}
According to \cite{Liang1999}, under certain assumptions, the parameter estimators possess both consistency and asymptotic normality.

\begin{lemma}\citep{Liang1999}\label{lem1}
	Suppose that there exists an $s > 2$ such that $E\|\mathbf{X}_{i}\|^{2s} < \infty$, $E\|\mathbf{U}_{i}\|^{2s} < \infty$ and $E|\epsilon_{i}|^{2s} < \infty$, where $\|\cdot\|$ represents an Euclidean norm.
	Additionally, {suppose that} $\frac{1}{n}\sum_{i=1}^{n}\mathbf{X}_{i}\mathbf{X}_{i}^{T}$ converges to the covariance matrix $\mathcal{H}$, where $\mathcal{H}$ is a non-random positive definite matrix. Then it follows that\\
	$(\mathrm{i})$ $\|\hat{\boldsymbol{\beta}} - \boldsymbol{\beta}\| = O_{P}(n^{-\frac{1}{2}})$ as $n\rightarrow \infty$.\\
	$(\mathrm{ii})$ $\hat{\boldsymbol{\beta}}$ follows an asymptotic normal distribution, i.e.,
	$$
	\sqrt{n}(\hat{\boldsymbol{\beta}} - \boldsymbol{\beta})\xrightarrow{d} N_{p}(\mathbf{0}, \mathcal{H}^{-1}\Gamma \mathcal{H}^{-1}),
	$$
	as $n\rightarrow \infty$, where
	$
	\Gamma =\sigma^{2}\mathcal{H} + \mathcal{H}\boldsymbol{\beta}^{T}\boldsymbol{\Sigma}_{uu}\boldsymbol{\beta} + E[(\mathbf{U}\mathbf{U}^{T} - \boldsymbol{\Sigma}_{uu})\boldsymbol{\beta}]^{\otimes 2} + \sigma^{2}\boldsymbol{\Sigma}_{uu}$,
the expression $O_{P}(1)$ denotes a sequence that is bounded in probability, $\xrightarrow{d}$ denotes convergence in distribution, $N_{p}$ denotes $p$-dimensional multivariate normal distribution,	and $\boldsymbol{v}^{\otimes 2} = \boldsymbol{vv}^{T}$ for any vector $\boldsymbol{v}$.
\end{lemma}

\begin{remark}\label{re1}
	In applications, $\boldsymbol{\Sigma}_{uu}$ is often unknown. \cite{Carroll2006, LiangLi2009} proposed that the estimate $\hat{\boldsymbol{\Sigma}}_{uu}$ of $\boldsymbol{\Sigma}_{uu}$ can be obtained by utilizing repeated observations corresponding to each $\mathbf{X}_{i}$. Let $J_{i}$ is times of repetition.
	For repetitive observations $\mathbf{W}_{i,j} = \mathbf{X}_{i} + \mathbf{U}_{i,j}, j = 1, \ldots, J_{i}, i = 1, \ldots, n $, we have
	$$\widehat{\boldsymbol{\Sigma}}_{uu} = \frac{\sum_{i=1}^{n} \sum_{j=1}^{J_{i}}(\mathbf{W}_{i,j} - \overline{\mathbf{W}}_{i}) (\mathbf{W}_{i,j} - \overline{\mathbf{W}}_{i})^{T}}{\sum_{i=1}^{n}(J_{i} - 1)}, $$
	where $\overline{\mathbf{W}}_{i} = J_{i}^{-1}\sum_{j=1}^{J_{i}}\mathbf{W}_{i,j}$. By substituting $\mathbf{W}_{i}$ by $\overline{\mathbf{W}}_{i}$, the estimator $$\left\{\sum_{i=1}^{n} \left(\overline{\mathbf{W}}_{i}\overline{\mathbf{W}}_{i}^{T} -J_{i}^{-1}\widehat{\boldsymbol{\Sigma}}_{uu}\right)\right\}^{-1} \sum_{i=1}^{n}\overline{\mathbf{W}}_{i}y_{i}$$ is unbiased and consistent, and the aforementioned properties in Lemma \ref{lem1} remain valid. We will use it in Section 6.2.
\end{remark}

\section{Optimal Subsampling Based on Corrected Likelihood}

In this section, we present a general subsampling algorithm that can be used to obtain estimators with measurement error models. Then we establish the consistency and asymptotic normality of these estimators and calculate the optimal subsampling probability using the A (or L)-optimality criterion.

\subsection{General subsampling algorithm}

Let us introduce some symbols before proposing the algorithm. We take a random subsample from the full data with replacement based on the subsampling probability $\pi_{i}, i=1, \ldots, n$, where $\sum_{i=1}^{n}\pi_{i} = 1$. In addition, we denote the subsample as $\mathcal{F}_{r} = \{(\mathbf{W}_{i}^{*}, y_{i}^{*})\}_{i=1}^{r}$ with $r$ being the subsample size. The corresponding subsampling probability is denoted as $\pi_{i}^{*}, i = 1, \ldots, r$. 
We can formulate the weighted loss function as follows
$$\ell^{*}(\boldsymbol{\beta})= \frac{1}{2n}\sum_{i=1}^{r} \frac{1}{r \pi_{i}^{*}}( y_{i}^{*} - \mathbf{W}_{i}^{*T} \boldsymbol{\beta} )^{2} -  \frac{1}{2}\boldsymbol{\beta}^{T} \boldsymbol{\Sigma}_{uu} \boldsymbol{\beta}.$$

\begin{algorithm}[h!]
	\caption{ General subsampling algorithm based on corrected likelihood}
	\label{alg1}
	\begin{itemize}
		\item {\bf 1. Subsample:} Extract a random subsample $\mathcal{F}_{r} = \{(\mathbf{W}_{i}^{*}, y_{i}^{*})\}_{i=1}^{r}$ of size $r$ ($r\ll n$) from the full dataset with replacement based on the subsampling probabilities $\{\pi_{i}\}_{i=1}^{n}$. The subsampling probabilities corresponding to the subsample are denoted as $\{\pi_{i}^{*}\}_{i=1}^{r}$. \par
		\item {\bf 2. Estimate:} Minimize the objective function $\ell^{*}(\boldsymbol{\beta})$ based on the subsample $\mathcal{F}_{r}$ to obtain the parameter estimator, that is,
		$\tilde{\boldsymbol{\beta}}  = \arg\underset{\boldsymbol{\beta}}{\min} \ell^{*}(\boldsymbol{\beta}).$
	\end{itemize}
\end{algorithm}

Algorithm \ref{alg1} is a general subsampling algorithm that relies on the corrected likelihood to address big data subsampling problems with measurement errors.

\subsection{Asymptotic properties}

In order to obtain the asymptotic properties of $\tilde{\boldsymbol{\beta}}$, the following assumptions are {needed}.

\begin{assumption} \label{js1}
	The matrix $\mathcal{H}_{W}=\frac{1}{n}\sum_{i=1}^{n}\mathbf{W}_{i}\mathbf{W}_{i}^{T}-\boldsymbol{\Sigma}_{uu}$ converges to a positive definite matrix in probability.
\end{assumption}

\begin{assumption} \label{js2}
	$\frac{1}{n} \sum_{i=1}^{n} \|\mathbf{W}_{i}\|^{4} = O_{P}(1)$ and $\frac{1}{n} \sum_{i=1}^{n} (y_{i} - \mathbf{W}_{i}^{T}\hat{\boldsymbol{\beta}})^{4} = O_{P}(1)$, where $\|\mathbf{W}_{i}\|^{k} = (\sum_{j=1}^{p}w_{ij}^{2})^{k/2}$.
\end{assumption}

\begin{assumption} \label{js3}
	$\max\limits_{i=1,\ldots,n} (n\pi_{i})^{-1} = O_{P}(1).$
\end{assumption}

\begin{assumption} \label{js4}
	There exists {some} $\delta > 0$ such that
	$\frac{1}{n} \sum_{i=1}^{n} (y_{i} - \mathbf{W}_{i}^{T} \hat{\boldsymbol{\beta}})^{2+\delta} \| \mathbf{W}_{i} \|^{2+\delta} = O_{P}(1).$
\end{assumption}

Assumption \ref{js1} indicates that $E(\mathbf{X}\mathbf{X}^{T})$ is positive definite and consistent with \cite{WangZhuMa2018}. Assumption \ref{js2} requires the limited moment, and ensures the consistency of parameter estimator. It is a special case of assumptions (H.3) and (H.6) in \cite{AiYuZhangWang2021}.
Assumption \ref{js3} restricts the weights in $\ell^{*}(\boldsymbol{\beta})$ in order to protect the estimating equation from being dominated by data points with minimal subsampling probabilities.
Assumption \ref{js4} introduces a finite moment, which ensures the asymptotic normality of the parameter estimator. It is also a special case of assumption (H.7) in \cite{AiYuZhangWang2021}.

\begin{theorem}\label{theo1}
	Under Assumptions \ref{js1}--\ref{js3}, as $r\rightarrow\infty$ and $n\rightarrow\infty$, $\tilde{\boldsymbol{\beta}}$ converges to $\hat{\boldsymbol{\beta}}$ in conditional probability given $\mathcal{F}_{n}$ and the convergence rate is $\sqrt{r}$. That is, with probability approaching one, for any $\varepsilon > 0$, there exist {finite} constants $\Delta_{\varepsilon}$ and $r_{\varepsilon}$ such that
	$$P(\|\tilde{\boldsymbol{\beta}} - \hat{\boldsymbol{\beta}}\| \geq r^{-1/2} \Delta_{\varepsilon} | \mathcal{F}_{n}) < \varepsilon$$
	for all $r> r_{\varepsilon}$.
\end{theorem}

\begin{remark}
	If a sequence of random variables is bounded in conditional probability, it will also be bounded in unconditional probability. Therefore, Theorem \ref{theo1} implies that $\|\tilde{\boldsymbol{\beta}} - \hat{\boldsymbol{\beta}}\| = O_{P}(r^{-1/2})$ (see \cite{WangZhuMa2018}).
\end{remark}

\begin{theorem}\label{theo2}
	Under Assumptions \ref{js1}--\ref{js4}, if $r = o(n)$, then as $r \rightarrow\infty$ and $n \rightarrow\infty$, $\tilde{\boldsymbol{\beta}} - \hat{\boldsymbol{\beta}}$ converges to a normal distribution in conditional probability given $\mathcal{F}_{n}$, that is,
	$$V^{-1/2} (\tilde{\boldsymbol{\beta}} - \hat{\boldsymbol{\beta}}) \xrightarrow{d} N_{p}(\mathbf{0},I),$$
	where the asymptotic covariance matrix
	$V = \mathcal{H}_{W}^{-1} V_{c} \mathcal{H}_{W}^{-1}$,
	and
	$$V_{c} = \frac{1}{r n^{2}} \sum_{i=1}^{n} \frac{(y_{i} - \mathbf{W}_{i}^{T}\hat{\boldsymbol{\beta}})^{2}\mathbf{W}_{i} \mathbf{W}_{i}^{T}}{\pi_{i}} -  \frac{1}{r}(\boldsymbol{\Sigma}_{uu}\hat{\boldsymbol{\beta}})^{\otimes 2},$$
	where $\boldsymbol{v}^{\otimes 2} = \boldsymbol{vv}^{T}$ for any vector $\boldsymbol{v}$.
\end{theorem}

\begin{remark}
	Theorem \ref{theo2} shows that the second term of $V_{c}$ is independent of the sampling probabilities when $\mathcal{F}_{n}$ is given. Hence it can be disregarded when we compute the optimal subsampling probabilities.
\end{remark}

\subsection{Optimal subsampling algorithm}

To determine the optimal subsampling probabilities, we utilize the A (or L)-optimality criterion from optimal design of experiments. This criterion aims to minimize the asymptotic mean {squared} error of $\tilde{\boldsymbol{\beta}}$ (or $\mathcal{H}_{W}\tilde{\boldsymbol{\beta}}$ ). Given that $\tilde{\boldsymbol{\beta}}$ (or $\mathcal{H}_{W}\tilde{\boldsymbol{\beta}}$ ) is asymptotically unbiased, minimizing the asymptotic variance $V$ (or $V_{c}$) is sufficient.

\begin{theorem}\label{theo3}
In Algorithm \ref{alg1}, we can draw two conclusions as follows.

(i) In Algorithm \ref{alg1}, if the subsampling probabilities $\pi_{i}$ for $i= 1,\ldots , n$ are selected as
		\begin{equation}\label{e10}
			\pi_{i}^{mV} = \frac{|y_{i} - \mathbf{W}_{i}^{T} \hat{\boldsymbol{\beta}}| \|\mathcal{H}_{W}^{-1} \mathbf{W}_{i}\|}{\sum_{i=1}^{n} |y_{i} - \mathbf{W}_{i}^{T} \hat{\boldsymbol{\beta}}| \|\mathcal{H}_{W}^{-1} \mathbf{W}_{i}\|}, i = 1, \ldots, n,
		\end{equation}
		then the asymptotic variance $tr(V)$ of $\tilde{\boldsymbol{\beta}}$ attains its minimum, which corresponds to the A-optimality criterion.
		
(ii)	In Algorithm \ref{alg1}, if the subsampling probabilities $\pi_{i}$ for $i= 1,\ldots , n$ are selected as
		\begin{equation}\label{e11}
			\pi_{i}^{mV_{c}} = \frac{|y_{i} - \mathbf{W}_{i}^{T} \hat{\boldsymbol{\beta}}| \|\mathbf{W}_{i}\|}{\sum_{i=1}^{n} |y_{i} - \mathbf{W}_{i}^{T} \hat{\boldsymbol{\beta}}| \|\mathbf{W}_{i}\|}, i = 1, \ldots, n,
		\end{equation}
		then the asymptotic variance $tr(V_{c})$ of $\mathcal{H}_{W}\tilde{\boldsymbol{\beta}}$ attains its minimum, which corresponds to the L-optimality criterion.
\end{theorem}

In (\ref{e10}) and (\ref{e11}), we observe that the optimal subsampling probabilities with measurement errors are similar to those of the generalized linear model developed by {\cite{AiYuZhangWang2021}}. However, our A-optimal subsampling probability includes a correction term $\boldsymbol{\Sigma}_{uu}$ for $\mathcal{H}_{W}$. Notably, the computing time required to determine $\|\mathcal{H}_{W}^{-1}\mathbf{W}_{i}\|$ has a complexity of $O(np^{2})$. While the L-optimal subsampling method demands only $O(np)$ time to calculate $\|\mathbf{W}_{i}\|$. Due to the presence of $\hat{\boldsymbol{\beta}}$ in the optimal sampling probability, we adopt the two-step algorithm, which is summarized in Algorithm \ref{alg2}.

\begin{algorithm}[h!]
	\caption{Optimal subsampling algorithm based on corrected likelihood}
	\label{alg2}
	\begin{itemize}
		\item {\bf 1.} Apply the uniform subsampling probability ${1}/{n}$ to Algorithm \ref{alg1} with subsample size $r_{0}$ to obtain a pilot estimate $\tilde{\boldsymbol{\beta}}_{0}$ of $\boldsymbol{\beta}$.
		Replace $\hat{\boldsymbol{\beta}}$ in $\{\pi^{mV}_{i}\}_{i=1}^{n}$ (or $\{\pi^{mV_{c}}_{i}\}_{i=1}^{n}$) with $\tilde{\boldsymbol{\beta}}_{0}$ to obtain the optimal subsampling probabilities $\pi_{i}(\tilde{\boldsymbol{\beta}})$. \par
		\item {\bf 2.} Take the optimal subsampling probabilities $\pi_{i}(\tilde{\boldsymbol{\beta}})$ into Algorithm \ref{alg1} with subsample size $r$. Combine the $r_{0} + r$ samples obtained from the two steps, and obtain $\breve{\boldsymbol{\beta}}$.
	\end{itemize}
\end{algorithm}

\begin{theorem}\label{theo4}
	Under Assumptions \ref{js1}--\ref{js3}, as $r_{0}r^{-1/2}\rightarrow 0, r\rightarrow\infty$ and $n\rightarrow\infty$, if $\tilde{\boldsymbol{\beta}}_{0}$ exists, then $\breve{\boldsymbol{\beta}}$ converges to $\hat{\boldsymbol{\beta}}$ in conditional probability given $\mathcal{F}_{n}$, and the convergence rate is $\sqrt{r}$. That is, with probability approaching one, for any $\varepsilon > 0$, there exist constants {finite} $\Delta_{\varepsilon}$ and $r_{\varepsilon}$ such that
	$$P(\|\breve{\boldsymbol{\beta}} - \hat{\boldsymbol{\beta}}\| \geq r^{-1/2} \Delta_{\varepsilon} | \mathcal{F}_{n}) < \varepsilon $$
	for all $r> r_{\varepsilon}$, where $\breve{\boldsymbol{\beta}}$ is obtained by Algorithm \ref{alg2}.
\end{theorem}

\begin{theorem}\label{theo5}
	Under Assumptions \ref{js1}--\ref{js4}, if $r = o(n)$, then as $r_{0}r^{-1/2}\rightarrow 0$ and $r \rightarrow\infty$ , $\breve{\boldsymbol{\beta}} - \hat{\boldsymbol{\beta}}$ converges to a normal distribution in conditional probability given $\mathcal{F}_{n}$, that is,
	$$(V^{\tilde{\boldsymbol{\beta}}_{0}})^{-1/2} (\breve{\boldsymbol{\beta}} - \hat{\boldsymbol{\beta}}) \xrightarrow{d} N_{p}(\mathbf{0},I),$$
	where $\breve{\boldsymbol{\beta}}$ is obtained by Algorithm \ref{alg2},
	$V^{\tilde{\boldsymbol{\beta}}_{0}} = \mathcal{H}_{W}^{-1} V_{c}^{\tilde{\boldsymbol{\beta}}_{0}} \mathcal{H}_{W}^{-1},$
	and
	$$V_{c}^{\tilde{\boldsymbol{\beta}}_{0}} = \frac{1}{r n^{2}} \sum_{i=1}^{n} \frac{(y_{i} - \mathbf{W}_{i}^{T}\hat{\boldsymbol{\beta}})^{2}\mathbf{W}_{i} \mathbf{W}_{i}^{T}}{\pi_{i}(\tilde{\boldsymbol{\beta}}_{0})} - \frac{1}{r}(\boldsymbol{\Sigma}_{uu}\hat{\boldsymbol{\beta}})^{\otimes 2}.$$
\end{theorem}

The standard error of an estimator has significant importance for statistical inference, such as hypothesis test and constructing confidence intervals. However, the computation of asymptotic covariance matrices necessitates utilizing the full data, which can be highly resource-intensive given the substantial sample size. To reduce computational costs, we employ the subsample to approximate the covariance matrix of $\breve{\boldsymbol{\beta}}$. This approximation is denoted as
$$\breve{V} = \breve{\mathcal{H}}_{W}^{-1} \breve{V}_{c} \breve{\mathcal{H}}_{W}^{-1},$$
where
$\breve{\mathcal{H}}_{W} = \frac{1}{n(r_{0}+r)} \sum_{i=1}^{r+r_{0}} \frac{\mathbf{W}_{i}^{*}\mathbf{W}_{i}^{*T}}{\pi_{i}^{*}} - \boldsymbol{\Sigma}_{uu},$
and
\begin{align*}
	\breve{V}_{c} =& \frac{1}{(r+r_{0})^{2} n^{2}} \sum_{i=1}^{r+r_{0}} \frac{(y_{i}^{*} - \mathbf{W}_{i}^{*T}\breve{\boldsymbol{\beta}})^{2}\mathbf{W}_{i}^{*} \mathbf{W}_{i}^{*T}}{(\pi_{i}^{*})^{2}} \\ &- \frac{1}{r+r_{0}}(\boldsymbol{\Sigma}_{uu}\breve{\boldsymbol{\beta}})^{\otimes 2}.
\end{align*}
Here, $\breve{\mathcal{H}}_{W}$ and $\breve{V}_{c}$ represent moment estimations for $\mathcal{H}_{W}$ and $V_{c}$, respectively. If $\breve{\boldsymbol{\beta}}$ is substituted with $\hat{\boldsymbol{\beta}}$, these estimates become unbiased.

\section{Perturbation Subsampling Based on Corrected Likelihood}

The optimal subsampling algorithm mentioned in Section 3 requires the calculation of unequal sampling probabilities for the full data {simultaneously}. However, as the sample size increases, implementation becomes increasingly memory-intensive. Therefore, this section {proposes} the perturbation subsampling algorithm to address this issue.

\subsection{Perturbation subsampling algorithm}

Suppose that
$\{\mu_{i}\}_{i=1}^{n}$ are generated from a Bernoulli distribution with probability $q_{n}=r/n$, $\{\nu_{i}\}_{i=1}^{n}$ are generated from a known probability distribution with mean $1/q_{n}$ and { the corresponding variance $b_{n}^{2}>0$}. The weighted loss function can be written as
$$L^{*}(\boldsymbol{\beta}) = \frac{1}{2n}\sum_{i=1}^{n} \psi_{i} (y_{i} - \mathbf{W}_{i}^{T} \boldsymbol{\beta} )^{2} - \frac{1}{2} \boldsymbol{\beta}^{T} \boldsymbol{\Sigma}_{uu} \boldsymbol{\beta},$$
where $\psi_{i} = \mu_{i}\nu_{i}$.

\begin{algorithm}[h!]
	\caption{ Perturbation subsampling algorithm based on corrected likelihood}
	\label{alg3}
	\begin{itemize}
		\item Set $k=1$.
		\item {\bf 1. Sampling:} Generate $n$ i.i.d. random variables $\{\mu_{k,i}\}_{i=1}^{n}$ such that $\mu_{k,i} \sim Bernoulli(q_{n})$, where $q_{n} = r/n$.
		\item {\bf 2. Random weighting:} Generate $n$ i.i.d. random variables $\{\nu_{k,i}\}_{i=1}^{n}$ from a completely known distribution with $E(\nu_{k,i}) = 1/q_{n}$ and variance being $b_{n}^{2}$.
		\item {\bf 3. Estimation:} Minimize $L^{*}(\boldsymbol{\beta})$ to obtain parameter estimator $\check{\boldsymbol{\beta}}_{k} = \arg\underset{\boldsymbol{\beta}}{\min} L^{*}(\boldsymbol{\beta}).$
		\item {\bf 4. Combination:} {Set $k=k+1$, and repeat Step 1--3 until $k=m$,} then combine the resulting estimates, namely $\check{\boldsymbol{\beta}}^{(m)} = \frac{1}{m} \sum_{k=1}^{m} \check{\boldsymbol{\beta}}_{k}.$
	\end{itemize}
\end{algorithm}

A general repeatedly random perturbation subsampling algorithm is given in Algorithm \ref{alg3}. The conditional variance of $\check{\boldsymbol{\beta}}^{(m)}$ can be estimated as
$\widehat{Var}(\check{\boldsymbol{\beta}}^{(m)}|\mathcal{F}_{n}) = \frac{1}{m(m-1)}\sum_{k=1}^{m}\left(\check{\boldsymbol{\beta}}_{k} - \check{\boldsymbol{\beta}}^{(m)}\right)\left(\check{\boldsymbol{\beta}}_{k} - \check{\boldsymbol{\beta}}^{(m)}\right)^{T}.$

\subsection{Asymptotic properties}

To obtain the asymptotic properties of $\check{\boldsymbol{\beta}}^{(m)}$, the following assumption is {necessary}.

\begin{assumption}\label{js6}
	$\limsup\limits_{n\rightarrow \infty} q_{n}E(\psi^{2}) < \infty$, and there {exists some} $\alpha>0$ such that $\limsup\limits_{n\rightarrow \infty}q_{n}^{2 + \alpha}E\nu^{2+ \alpha}<\infty$.
\end{assumption}

In Assumption \ref{js6}, finite second-order moment and higher-order moment are assumed, {it} is equivalent to that there exists some $\alpha>0$ such that $\limsup\limits_{n\rightarrow \infty} q_{n}^{1 + \alpha} E \psi^{2+ \alpha} < \infty$.

\begin{theorem}\label{theo6}
	Under Assumptions \ref{js1}, \ref{js2} {and} \ref{js6}, as $r \rightarrow \infty $ and $ n \rightarrow \infty $, then $\check{\boldsymbol{\beta}}^{(m)}$ converges to $\hat{\boldsymbol{\beta}}$ in conditional probability given $\mathcal{F}_{n}$, and the convergence rate is $(mr)^{1/2}$. That is, with probability approaching one, for any $\varepsilon > 0$, there exist constants $\Delta_{\varepsilon}$ and $r_{\varepsilon}$ such that
	$$
	P(\|\check{\boldsymbol{\beta}}^{(m)} - \hat{\boldsymbol{\beta}}\| \geq (mr)^{-1/2}\Delta_{\varepsilon} | \mathcal{F}_{n}) < \varepsilon
	$$
	for all $r> r_{\varepsilon}$.
\end{theorem}

\begin{remark}
	According to Theorem \ref{theo6}, $\check{\boldsymbol{\beta}}^{(m)}$ is the consistent estimator of $\hat{\boldsymbol{\beta}}$. When $rm < n$, the convergence rate is $(rm)^{1/2}$, otherwise it is $n^{1/2}$. Hence, the estimation based on the full data is still more effective than that using repeat perturbation subsampling. Intuitively, the convergence rate of $\check{\boldsymbol{\beta}}^{(m)} (m>1)$ is faster than $\check{\boldsymbol{\beta}}^{(1)}$. So the repeated perturbation subsampling method produces smaller estimation errors for the same subsample size compared to the optimal subsampling algorithm. Note that the stochastic weights $\psi_{i}$ in Algorithm \ref{alg3} are independent of the data, so this algorithm can be implemented under parallel or
distributed computing framework. Simulation results declare that it takes
less computational time than the optimal subsampling algorithm.
\end{remark}

\begin{theorem}\label{theo7}
	Under Assumptions \ref{js1}, \ref{js2}, \ref{js4} {and} \ref{js6}, if $r = o(n)$, then as $r \rightarrow\infty$ and $n \rightarrow\infty$, $\check{\boldsymbol{\beta}}^{(m)} - \hat{\boldsymbol{\beta}}$ converges to a normal distribution in conditional probability given $\mathcal{F}_{n}$, that is, (i) when $mr<n$,
	$$ \Sigma^{-1/2} \sqrt{m r/a_{n}}(\check{\boldsymbol{\beta}}^{(m)} - \hat{\boldsymbol{\beta}}) \xrightarrow{d} N_{p}(\mathbf{0},I), r \rightarrow \infty , n \rightarrow \infty;$$ (ii) when $mr\geq n$,
	$$ \Sigma^{-1/2} \sqrt{n/a_{n}}(\check{\boldsymbol{\beta}}^{(m)} - \hat{\boldsymbol{\beta}}) \xrightarrow{d} N_{p}(\mathbf{0},I), r \rightarrow \infty , n \rightarrow \infty ,$$
	where $a_{n} = 1-q_{n} +b_{n}^{2} q_{n}^{2}$, $\Sigma = \mathcal{H}_{W}^{-1} \Sigma_{c} \mathcal{H}_{W}^{-1},$
	and
	$\Sigma_{c} = \frac{1}{n} \sum_{i=1}^{n} \mathbf{W}_{i}\mathbf{W}_{i}^{T} (y_{i} - \mathbf{W}_{i}^{T}\hat{\boldsymbol{\beta}})^{2}.$
\end{theorem}

\section{Simulation Studies}


We generate the full data from model (\ref{eq1}) with $n=10^{4}$, $\boldsymbol{\beta}=(1,1,1,1,1)^T$
and $\epsilon_{i} \sim N(0,1)$.
We consider the following two cases to generate the covariates $\mathbf{X}_i$ \citep{WangZhuMa2018}.

Case 1: $\mathbf{X}_i \sim N_5(\mathbf{0},\boldsymbol{\Sigma})$ {with} $\boldsymbol{\Sigma}_{j,k} = 0.5^{|j-k|}$ for $j, k=1, 2, \ldots, 5$.

Case 2: $\mathbf{X}_i \sim t_{3}(\mathbf{0},\boldsymbol{\Sigma})$ {with} $\boldsymbol{\Sigma}_{j,k} = 0.5^{|j-k|}$ for $j, k=1, 2, \ldots, 5$.

Let $\mathbf{U}_i \sim N_5(\mathbf{0},\sigma_{u}^{2}I)$, then $\mathbf{W}_{i} = \mathbf{X}_{i} + \mathbf{U}_{i}$.
We consider the following three values for $\sigma_{u}^{2}$ respectively:
$\sigma_{u}^{2} = 0.6, 0.4, 0.2$;

For Algorithm \ref{alg3}, unless otherwise specified, fix $m = 10$, and assume that the known distribution with random weighting follows an exponential distribution with mean $1/q_{n}$, i.e., $\nu_{i} \sim \mathrm{Exp}(q_{n})$. Correspondingly, $b_{n}^{2} = 1/q_{n}^{2}$ and $a_{n} = 1- q_{n} +b_{n}^{2}q_{n}^{2} = 2 - q_{n}$.
We choose $r_{0}=500$ and $r=500, 1000, 1500, 2000, 3000$. For each value of $r$, we perform $N=1000$ repetitions to calculate the mean squared error (MSE): $\frac{1}{N} \sum_{i=1}^{N}\|\hat{\boldsymbol{\beta}}_{i} - \boldsymbol{\beta}\|^{2}$.

\begin{figure*}[htpb]
	\centering{ \includegraphics[scale=0.45]{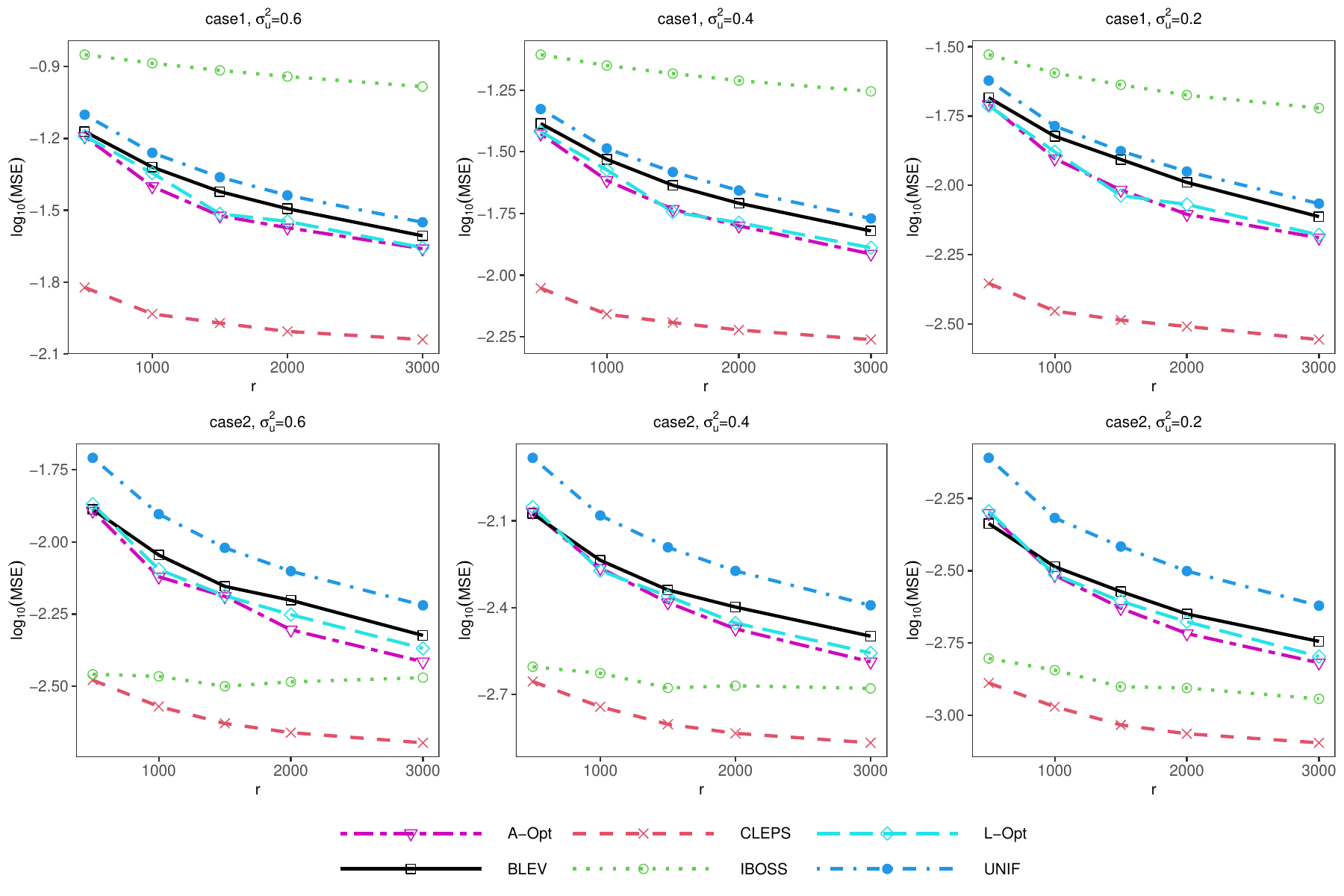} }
	\caption{The MSEs based on different $X$ and $\sigma_{u}^{2}$ for different $r$.}
	\label{tu1}
\end{figure*}

\begin{figure*}[htpb]
	\centering{ \includegraphics[scale=0.45]{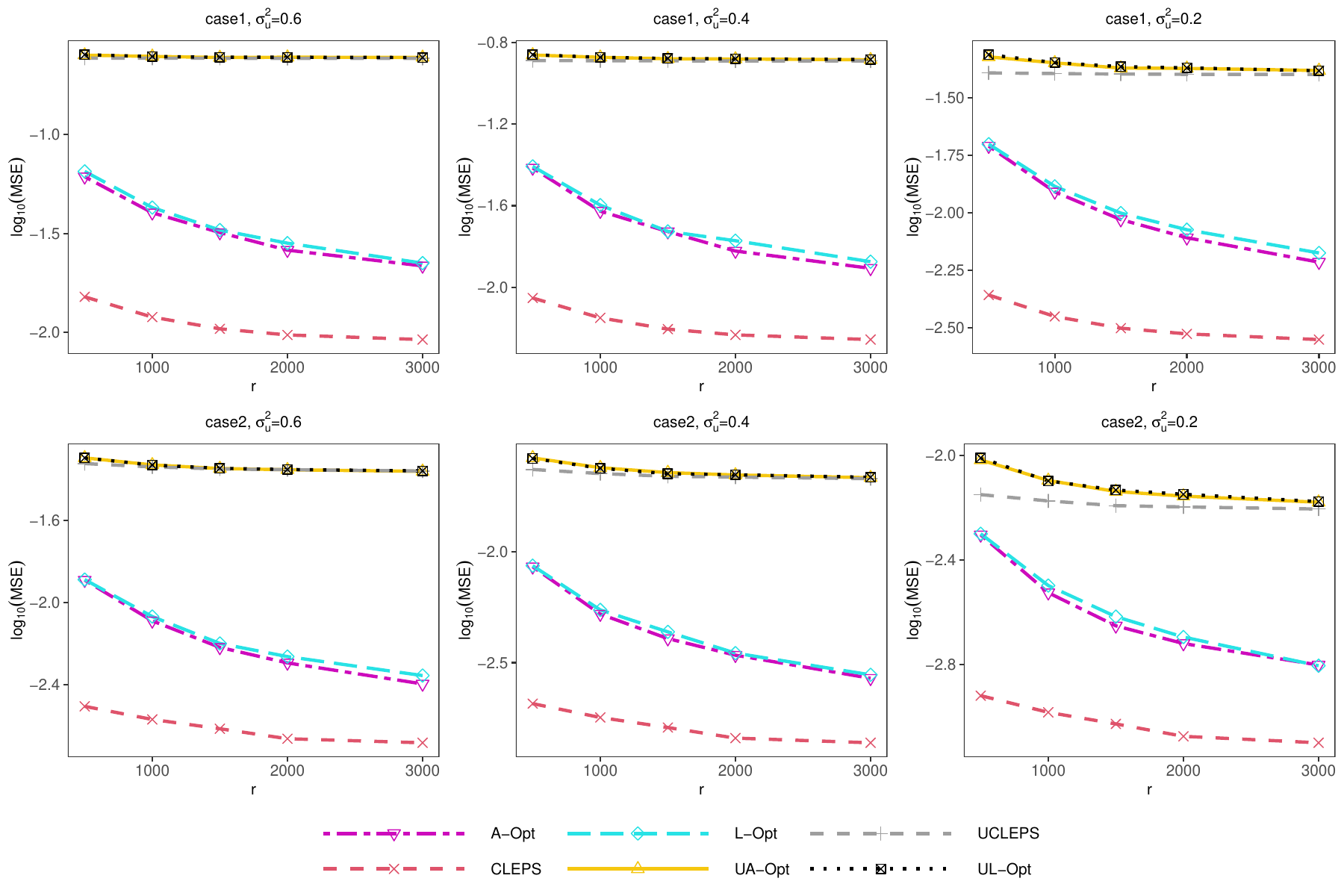} }
	\caption{The MSEs based on different $X$ and $\sigma_{u}^{2}$ for different $r$.}
	\label{tu2}
\end{figure*}

In Figure \ref{tu1}, we consider six corrected subsampling methods with measurement errors, including the perturbation subsampling based on corrected likelihood in Algorithm \ref{alg3} (CLEPS), A (or L)-optimal subsampling based on corrected likelihood in Algorithm \ref{alg2} (A-Opt and L-Opt), corrected uniform subsampling (UNIF), corrected leverage-based subsampling (BLEV), and corrected D-optimal subsampling (IBOSS). In Figure \ref{tu2}, we compare CLEPS, A-Opt, and L-Opt with the corresponding uncorrected subsampling methods with measurement errors, including the uncorrected perturbation subsampling (UCLEPS), uncorrected A (or L)-optimal subsampling (UA-Opt and UL-Opt).
For fair comparisons, except for the optimal subsampling methods, all the other methods use $r_{0} + r$ subsamples for parameter estimation.

The results presented in Figure \ref{tu1} demonstrate the superior performance of CLEPS compared to other methods, closely followed by the A-Opt and L-Opt. Furthermore, as the subsample size increases, the MSEs of our proposed methods approach to zero, which confirm the inconsistency of ordinary least squares estimator for linear models with measurement errors.
Additionally, decreasing $\sigma_{u}^{2}$ lead to reductions for MSEs obtained by all six methods. For the IBOSS method, the decrease in the MSE is not as evident, and its advantage is more significant if the tail of the covariate distribution is heavier. These results are consistent with \cite{WangYang2019}. However, the performance of CLEPS, A-Opt and L-Opt does not seem to be affected by the covariate distribution. It is seen from Figure \ref{tu2} that our proposed methods perform much better than the uncorrected subsampling methods. It indicates that the corrected subsampling methods are efficient.

\begin{figure*}[htpb]
	\centering
	\begin{minipage}{0.45\linewidth}
		\centering
		\includegraphics[scale=0.6]{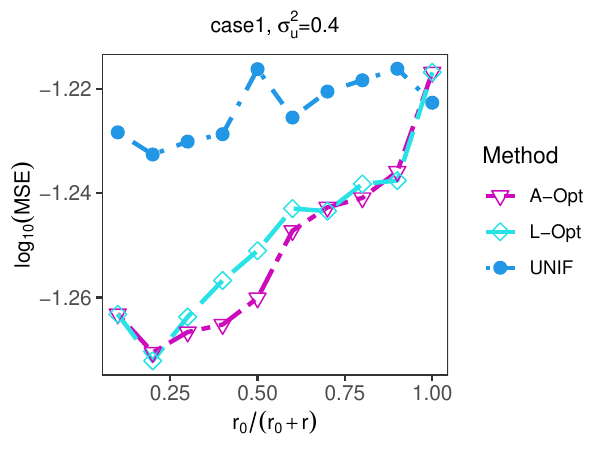}
		\caption*{(1)}
	\end{minipage}
	\begin{minipage}{0.45\linewidth}
		\centering
		\includegraphics[scale=0.6]{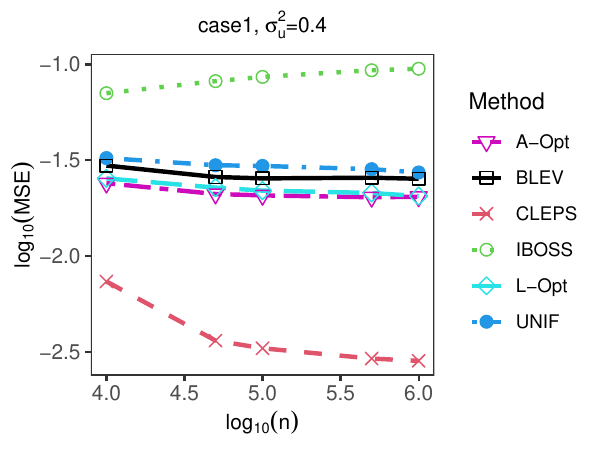}
		\caption*{(2)}
	\end{minipage}
	\begin{minipage}{0.45\linewidth}
		\centering
		\includegraphics[scale=0.6]{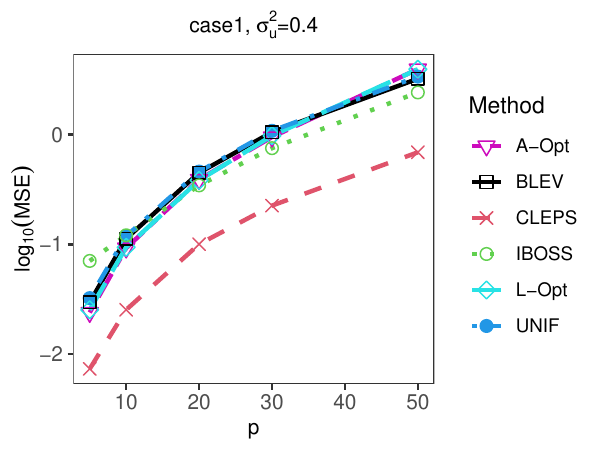}
		\caption*{(3)}
	\end{minipage}
	\begin{minipage}{0.45\linewidth}
		\centering
		\includegraphics[scale=0.6]{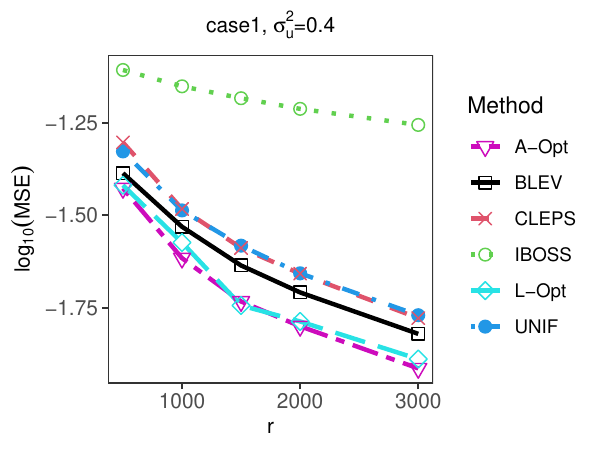}
		\caption*{(4)}
	\end{minipage}
	\caption{(1) The MSEs for different $r_{0}$ with fixed $r_{0}+r=1000$. (2) and (3) The MSEs for different $n$ and $p$ with $m=10, r_{0}=500, r=1000$, respectively. (4) The MSEs for different $r$ with $m=1$.}
	\label{tu8}
\end{figure*}

To further investigate the impact of other parameters on the sampling {methods}, Figure \ref{tu8} offers the variations in MSE across distinct values of $r_{0}, n, p,$ and $m$ while keeping $\sigma_{u}^{2} = 0.4$. Additionally, Table \ref{biao1} presents the results about different {values of} $m$ for CLEPS.

Figure \ref{tu8} (1) shows that as $r_{0}$ {increases}, the performance of A-Opt and L-Opt initially enhances and subsequently declines for fixed $r_{0} + r = 1000$. {The} reason is that the inaccurate estimation is obtained in the first step when $r_{0}$ is too small.
{Figure \ref{tu8} (2)} indicates that as the sample size increases, the CLEPS performs more efficiently, while {the MSEs of} other methods remain almost unchanged.
{Figure \ref{tu8} (3)} displays that as the dimension $p$ increases, the MSEs of various methods also increases. When $p>15$, A-Opt and L-Opt underperform compared to IBOSS, whereas the CLEPS consistently demonstrates the best performance.
{Figure \ref{tu8} (4)} shows that when $m=1$, the MSE of CLEPS is slightly larger than those of the BLEV, A-Opt and L-Opt, and is similar to that of UNIF.
\begin{table}[h!]
	\caption{The MSEs for different values of $m$.}
	\label{biao1}
	{\tabcolsep=3.5pt
		\begin{tabular}{cccccc}
			\toprule
			$\log_{10}\mathrm{MSE}$     &$m=1$     & $m=10$      & $m=20$       & $m=30$       & $m=50$    \\
			\midrule
			$r=500$ & $-$1.32  & $-$2.06    & $-$2.20     & $-$2.23     & $-$2.28   \\
			$r=1000$ & $-$1.49  & $-$2.14    & $-$2.25     & $-$2.29     & $-$2.31   \\
			$r=1500$ & $-$1.58  & $-$2.19    & $-$2.27     & $-$2.30     & $-$2.32   \\
			$r=2000$ & $-$1.64  & $-$2.23    & $-$2.30     & $-$2.32     & $-$2.33    \\
			$r=3000$& $-$1.77  & $-$2.26    & $-$2.31     & $-$2.33     & $-$2.34    \\
			\bottomrule
		\end{tabular}
	}
\end{table}

In Table \ref{biao1}, it is observed that as $m$ increases, the MSE decreases. However, the rate of reduction in MSE progressively diminishes. This suggests that $m$ should be significantly smaller than $r$ in order to achieve effective inference. It is advisable to set $m < r/10$, in accordance with the findings in \cite{ShangCheng2017}, \cite{Wang2019b}, and \cite{WangMa2021}, which suggest that the number of partitions should be significantly smaller than the sample size within each data partition.

We use the Sys.time() function to calculate the computing times of different methods. Following \cite{WangYang2019}, we show the computing times (in seconds) for different combinations of the full data
size $n$ and the number of covariates $p$ for Case 1 with $r_0 = 500$ and $r_1=1000$. Computations are carried out on a desktop running Window 10 with an AMD E2 processor and 8GB memory. For comparison, the computing times for using
the full data (FULL) are also presented.
As shown in Table \ref{biao2}, the uniform subsampling algorithm requires the least
computing time. Note that $\hat{\boldsymbol{\beta}}$ can be directly obtained, so the computing time of the full data is much smaller than the BLEV, A-Opt and L-Opt. The CLEPS compares favorably to the other subsampling methods except for UNIF.

\begin{table}[h!]
	\caption{The computing times (in seconds) for different combinations of $n$ and $p$.}
	\label{biao2}
	{\tabcolsep=1.6pt
		\begin{tabular}{rrrrrrrrr}
			\toprule
$n$    &$p$	   &UNIF    &BLEV     & IBOSS     & A-Opt      &L-Opt      & CLEPS   & FULL   \\ \midrule
$10^4$ &$300$  &0.092&  2.679 &0.721 & 2.015 &1.457 &0.720 &0.599  \\
$10^5$ &       &0.093& 23.115 &5.975 &13.050 &7.024 &0.731 &5.971  \\
$10^6$ &       &0.094&282.264 &67.482&159.018&76.134&1.392 &61.164  \\ \midrule
$10^5$ &$50$   &0.024&  1.297 &1.171 & 0.925 &0.505 &0.154 &0.302  \\
       &$100$  &0.013&  4.812 &2.588 & 2.277 &1.406 &0.267 &1.005  \\
       &$300$  &0.111& 31.669 &6.685 &19.766 &9.672 &0.831 &9.745  \\
\bottomrule
\end{tabular}
}
\end{table}

\section{ Real Examples}

{\color{red}
To illustrate the effectiveness of our method in application, we give two real data sets including the diamond price dataset and the airline delay dataset. For the diamond price dataset without repeated measurements, there are no measurement errors in covariates. So we artificially add the measurement errors to the covariates with a known variance-covariance matrix of the measurement error vector. For the airline delay dataset, the variance-covariance matrix of the measurement error vector can be estimated by the repeated measurements.}

\subsection{Diamond price dataset}

The diamond price dataset is an integrated dataset containing prices and other characteristics of approximately 54,000 diamonds. This dataset can be found at \url{https://www.kaggle.com/datasets/shivam2503/diamonds}. Our aim is to explore the relationship between diamond prices ($y$) and three covariates: carat ($x_{1}$), depth ($x_{2}$) and table ($x_{3}$). All variables are standardized and the linear model is
$$y = x_{1}\beta_{1} + x_{2}\beta_{2} + x_{3}\beta_{3} + \epsilon.$$

\begin{figure*}[htpb]
	\centering{ \includegraphics[scale=0.45]{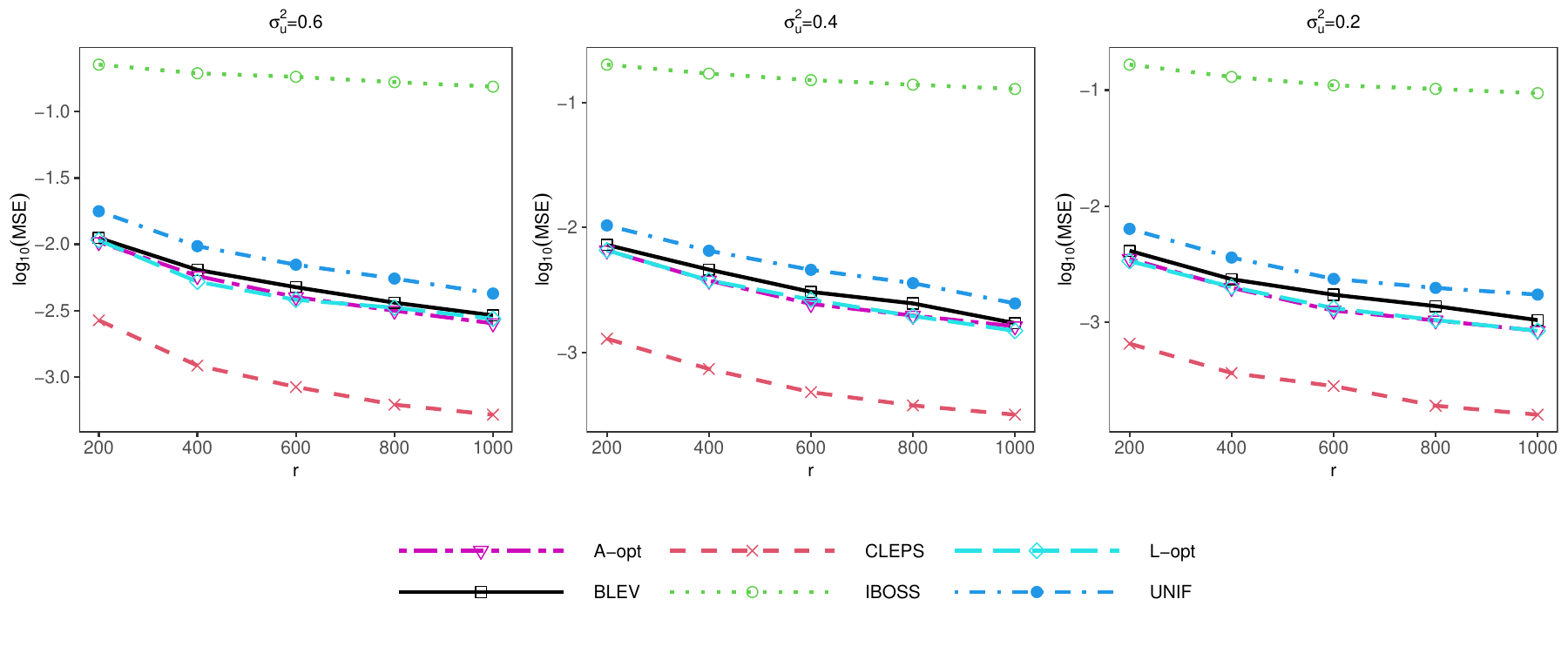} }
	\caption{The MSEs for different {values of} $\sigma_{u}^{2}$ and $r$ with $m=10$.}%
	\label{tu6}
\end{figure*}

Following \cite{Cui2025}, the measurement error model is $\mathbf{W} = \mathbf{X} + \mathbf{U}$, where $\mathbf{X}=(x_1,x_2,x_3)^T$, $\mathbf{U}$ is the measurement error vector with mean zero and covariance matrix $\sigma_{u}^{2}I_3$,  $\sigma_{u}^{2}= 0.6, 0.4, 0.2$.  Figure \ref{tu6} depicts MSEs for different {values of} $\sigma_{u}^{2}$ and $r$ when $m = 10, r_{0}=100$, {with repetition $N=500$}.
As the sample size {$r$} increases, the MSEs of the proposed methods decrease and become smaller than those of other methods.
Moreover, with a decrease in the variance of measurement error, there is a corresponding decrease in the MSEs of the proposed methods.

\subsection{Airline delay dataset}
The airline delay dataset has nearly 120 million records, which can be found {at} \url{https://community.amstat.org/jointscsg-section/dataexpo/dataexpo2009}. It includes detailed information on the arrival and departure of all commercial flights within the USA from 1987 to 2008. This study focuses solely on 2008 data, yielding a total of 2389217 samples and 29 variables. After cleaning, 715731 observation data points were obtained. We select Arrival Delay as the response variable, with Departure Delay, Distance, Air Time, and Elapsed Time as the covariates.  It's worth noting that the flight elapsed time incorporates two types of time: Actual Elapsed Time and CRS Elapsed Time. Among them, the CRS Elapsed Time refers to the original elapsed time. Therefore, the Flight Elapsed Time is regarded as a variable with measurement errors of two repetitions, while other variables have no measurement errors. {According to Remark \ref{re1}}, the estimate for the covariance matrix of measurement error is
\begin{equation*}
	\widehat{\boldsymbol{\Sigma}}_{uu} =
	\begin{bmatrix}
		0 &   & &  \\
		& 0 &  & \\
		&   & 0&\\
		&   &  &0.0195
	\end{bmatrix}.
\end{equation*}
In our proposed methods, the Flight Elapsed Time is calculated as the mean of the two variables, whereas in other methods, it is determined as the Actual Elapsed Time.

\begin{figure}[h!]
	\centering{ \includegraphics[scale=0.7]{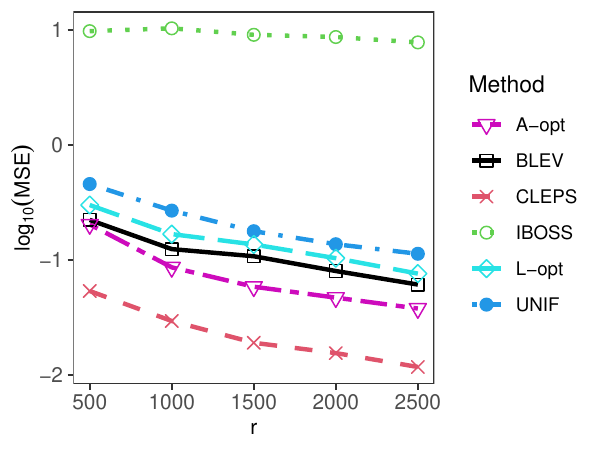} }
	\caption{The MSEs for different {values of $r$} with $m=10$.}%
	\label{tu7}
\end{figure}

Figure \ref{tu7} displays the impact of varying $r$ on MSE when $m = 10, r_{0} = 200$, and $N = 500$.
Notably, the MSEs of the proposed methods consistently decrease {with the increasing of $r$} and outperform other methods, indicating that our methods also have good properties in applications when $\boldsymbol{\Sigma}_{uu}$ is unknown.

\section{Conclusion and discussion}

The optimal subsampling methods have received much attention for big data analysis. However, much of the existing work on optimal subsampling has focused on the clean data case. In many applications, covariates may be observed inaccurately due to measurement errors. The numerical results show that the performance of the existing subsampling approaches is unsatisfactory when the covariates maybe inaccurately observed.
For large-scale data with measurement errors in covariates of linear models, this paper proposes two subsampling methods based on the corrected likelihood.
Theoretical results and numerical studies demonstrate that the two proposed algorithms outperform other existing sampling methods when measurement errors are present.
However, this work only considers linear models. Future research directions could include investigating nonlinear models, high-dimensional data, or distributed data.
\backmatter
\bmhead{Acknowledgements}
	The authors would like to thank three anonymous referees, and the Editor for their constructive comments that improved the
	quality of this paper. The authors are listed in alphabetical order.

\section*{Declarations}

\bmhead{Funding}
This research was supported by the National Natural Science Foundation of China
(12131001, 12171277, 12271294 and 12371260).

\bmhead{Author contributions}
Ju, J.: methodology, software, writing--original draft preparation;
Liu, M.: conceptualization, methodology, writing--review and editing;
Wang, M.: conceptualization, methodology, software, writing--review and editing, supervision;
Zhao, S.: conceptualization, methodology, writing--review and editing, supervision.

\bmhead{Conflict of interest}
The authors declare that there are no competing interests in this paper.

\begin{appendices}
\section{}
\begin{lemma}[$C_{r}$ inequality]\label{lem3}
	Let $\{X_{i}, i\geq 1\}$ be a sequence of independent random variables, then the $k$-th moment of the sum of random variables is not greater than the sum of the $k$-th moments of the random variables. i.e.,
	$$E\Big|\sum_{i=1}^{n}X_{i}\Big|^{k} \leq C_{r} \sum_{i=1}^{n} E|X_{i}|^{k}, $$
	where
	$
	C_{r}= \left\{
	\begin{aligned}
		&1,      & &0< k \leq 1,\\
		&n^{k-1},& &k > 1.
	\end{aligned}
	\right.
	$
\end{lemma}

\bmhead{Proof of Theorem 1}
\begin{proof}
	Note that
	\begin{equation}\label{e1}
		\begin{aligned}
			\tilde{\boldsymbol{\beta}} - \hat{\boldsymbol{\beta}}
			=& \left[\frac{1}{n} \sum_{i=1}^{r}\frac{1}{r\pi_{i}^{*}}\mathbf{W}_{i}^{*}\mathbf{W}_{i}^{*T} - \boldsymbol{\Sigma}_{uu}\right] ^{-1} \\
			&\cdot \left[\frac{1}{n}\sum_{i=1}^{r} \frac{1}{r\pi_{i}^{*}} \mathbf{W}_{i}^{*} (y_{i}^{*} - \mathbf{W}_{i}^{*T}\hat{\boldsymbol{\beta}}) + \boldsymbol{\Sigma}_{uu}\hat{\boldsymbol{\beta}}\right]\\
			=& -(\tilde{\mathcal{H}}_{W})^{-1} \dot{\ell}^{*}(\hat{\boldsymbol{\beta}}),
		\end{aligned}
	\end{equation}
	where
	$$\tilde{\mathcal{H}}_{W} = \frac{1}{r} \sum_{i=1}^{r}\left[\frac{1}{n\pi_{i}^{*}}\mathbf{W}_{i}^{*}\mathbf{W}_{i}^{*T} - \boldsymbol{\Sigma}_{uu}\right],$$
	$$\dot{\ell}^{*}(\boldsymbol{\beta}) = \frac{1}{n}\sum_{i=1}^{r} \frac{1}{r\pi_{i}^{*}} \left[-\mathbf{W}_{i}^{*} (y_{i}^{*} - \mathbf{W}_{i}^{*T}\boldsymbol{\beta})\right] - \boldsymbol{\Sigma}_{uu}\boldsymbol{\beta}. $$
	Therefore, we only need to prove
	\begin{equation}\label{e2}
		\dot{\ell}^{*}(\hat{\boldsymbol{\beta}}) = O_{P|\mathcal{F}_{n}}(r^{-1/2}),
	\end{equation}
	and
	\begin{equation}\label{e3}
		\tilde{\mathcal{H}}_{W} - \mathcal{H}_{W} = O_{P|\mathcal{F}_{n}}(r^{-1/2}),
	\end{equation}
	where $\mathcal{H}_{W} = \frac{1}{n} \sum_{i=1}^{n}\mathbf{W}_{i}\mathbf{W}_{i}^{T} - \boldsymbol{\Sigma}_{uu}. $
	
	To prove (\ref{e2}), we can directly obtain that
	\begin{align*}
		&E(\dot{\ell}^{*}(\hat{\boldsymbol{\beta}})|\mathcal{F}_{n})\\
		=& E\left\{\frac{1}{n}\sum_{i=1}^{r}\frac{1}{r \pi_{i}^{*}}\left[-\mathbf{W}_{i}^{*} (y_{i}^{*} - \mathbf{W}_{i}^{*T}\hat{\boldsymbol{\beta}})\right] - \boldsymbol{\Sigma}_{uu}\hat{\boldsymbol{\beta}} \bigg|\mathcal{F}_{n}\right\} \\
		=& \frac{1}{nr}\sum_{i=1}^{r} E\left\{\frac{1}{\pi_{i}^{*}}\left[-\mathbf{W}_{i}^{*} (y_{i}^{*} - \mathbf{W}_{i}^{*T}\hat{\boldsymbol{\beta}})\right]\Big|\mathcal{F}_{n}\right\}- \boldsymbol{\Sigma}_{uu}\hat{\boldsymbol{\beta}}\\
		=& \frac{1}{n}\sum_{i=1}^{n} \pi_{i}\cdot \frac{1}{\pi_{i}} \left[-\mathbf{W}_{i} (y_{i} - \mathbf{W}_{i}^{T}\hat{\boldsymbol{\beta}})\right] - \boldsymbol{\Sigma}_{uu}\hat{\boldsymbol{\beta}} \\
		=& \frac{1}{n}\sum_{i=1}^{n}\left[-\mathbf{W}_{i} (y_{i} - \mathbf{W}_{i}^{T}\hat{\boldsymbol{\beta}})\right] - \boldsymbol{\Sigma}_{uu}\hat{\boldsymbol{\beta}}\\
		=& \mathbf{0} .
	\end{align*}
	For the $j$-th element $\dot{\ell}_{j}^{*}(\hat{\boldsymbol{\beta}})$ of $\dot{\ell}^{*}(\hat{\boldsymbol{\beta}})$ {with} $1 \leq j \leq p $, 
	\begin{align*}
		&Var(\dot{\ell}_{j}^{*}(\hat{\boldsymbol{\beta}})|\mathcal{F}_{n})\\
		=&E\left\{\frac{1}{n}\sum_{i=1}^{r}\frac{1}{r \pi_{i}^{*}}\left[-w_{ij}^{*} (y_{i}^{*} - \mathbf{W}_{i}^{*T}\hat{\boldsymbol{\beta}})\right] - (\boldsymbol{\Sigma}_{uu}\hat{\boldsymbol{\beta}})_{j}\bigg|\mathcal{F}_{n}\right\}^{2}\\
		=&\frac{1}{r^{2}}\sum_{i=1}^{r}E\left\{\frac{1}{n \pi_{i}^{*}}\left[-w_{ij}^{*} (y_{i}^{*} - \mathbf{W}_{i}^{*T}\hat{\boldsymbol{\beta}})\right] - (\boldsymbol{\Sigma}_{uu}\hat{\boldsymbol{\beta}})_{j}\Big|\mathcal{F}_{n}\right\}^{2} \\
		=&\frac{1}{r} \Bigg\{ E\left\{\frac{1}{n \pi_{i}^{*}}\left[-w_{ij}^{*} (y_{i}^{*} - \mathbf{W}_{i}^{*T} \hat{\boldsymbol{\beta}})\right] \Big|\mathcal{F}_{n}\right\}^{2} \\
		&+ E\left\{(\boldsymbol{\Sigma}_{uu}\hat{\boldsymbol{\beta}})_{j}\Big|\mathcal{F}_{n}\right\}^{2}\\
		&-2E\left\{\frac{1}{n \pi_{i}^{*}}\left[-w_{ij}^{*} (y_{i}^{*} - \mathbf{W}_{i}^{*T}\hat{\boldsymbol{\beta}})\right] (\boldsymbol{\Sigma}_{uu}\hat{\boldsymbol{\beta}})_{j}\Big|\mathcal{F}_{n}\right\} \Bigg\}\\
		=&\frac{1}{n^{2}r}\sum_{i=1}^{n} \pi_{i}\cdot \left\{\frac{1}{\pi_{i}} \left[-w_{ij} (y_{i} - \mathbf{W}_{i}^{T}\hat{\boldsymbol{\beta}})\right]\right\}^{2} - \frac{1}{r}(\boldsymbol{\Sigma}_{uu}\hat{\boldsymbol{\beta}})_{j}^{2}\\
		=& \frac{1}{n^{2}r}\sum_{i=1}^{n} \frac{w_{ij}^{2} (y_{i} - \mathbf{W}_{i}^{T}\hat{\boldsymbol{\beta}})^{2}}{\pi_{i}}  - \frac{1}{r}(\boldsymbol{\Sigma}_{uu}\hat{\boldsymbol{\beta}})_{j}^{2}\\
		\leq & \frac{1}{r}\max_{i=1,\ldots,n}(n\pi_{i})^{-1}\sum_{i=1}^{n}\frac{\|\mathbf{W}_{i}\|^{2}(y_{i} - \mathbf{W}_{i}^{T}\hat{\boldsymbol{\beta}})^{2}}{n}\\
		&- \frac{1}{r}(\boldsymbol{\Sigma}_{uu}\hat{\boldsymbol{\beta}})_{j}^{2}.
	\end{align*}
	By Assumption \ref{js2} and the Holder inequality, we can obtain
		\begin{align}\label{ee1}
			&\sum_{i=1}^{n}\frac{\|\mathbf{W}_{i}\|^{2}(y_{i} - \mathbf{W}_{i}^{T}\hat{\boldsymbol{\beta}})^{2}}{n} \nonumber\\
			\leq& \left[\sum_{i=1}^{n}\frac{\|\mathbf{W}_{i}\|^{4}}{n}\right]^\frac{1}{2} \left[\sum_{i=1}^{n}\frac{(y_{i} - \mathbf{W}_{i}^{T}\hat{\boldsymbol{\beta}})^{4}}{n}\right]^\frac{1}{2}
			= O_{P}(1).
		\end{align}
	According to Assumption \ref{js3}, we can infer that
	$$Var(\dot{\ell}_{j}^{*}(\hat{\boldsymbol{\beta}})|\mathcal{F}_{n}) = \frac{1}{r}O_{P}(1)O_{P}(1)-O_{P}(r^{-1}) = O_{P}(r^{-1}).$$
	From the Chebyshev inequality, for a sufficiently large $M$, we have
	\begin{align*}
		P(\|\dot{\ell}^{*}(\hat{\boldsymbol{\beta}})\| \geq r^{-1/2}M |\mathcal{F}_{n})
		&\leq \frac{r E(\|\dot{\ell}^{*}(\hat{\boldsymbol{\beta}})\|^{2}|\mathcal{F}_{n})}{M^{2}} \\
		&= \frac{r \sum_{j=1}^{p}E(\dot{\ell}_{j}^{*}(\hat{\boldsymbol{\beta}})^{2}|\mathcal{F}_{n})}{M^{2}} \\
		&= \frac{O_{P}(1)}{M^{2}}\rightarrow 0 , n, r \rightarrow \infty.
	\end{align*}
	Thus, {equation} (\ref{e2}) is derived.
	
	In order to prove (\ref{e3}), we directly obtain that
	$$E(\tilde{\mathcal{H}}_{W} | \mathcal{F}_{n}) = \mathcal{H}_{W}.$$
	{For any element $\tilde{\mathcal{H}}_{W}^{j_{1}j_{2}}$ of $\tilde{\mathcal{H}}_{W}$, $1\leq j_{1}, j_{2} \leq p$, based on} Assumptions \ref{js2}, \ref{js3} and Cauchy-Schwarz inequality, it can be concluded that
	\begin{align*}
		&Var(\tilde{\mathcal{H}}_{W}^{j_{1}j_{2}} | \mathcal{F}_{n})\\
		=& E(\tilde{\mathcal{H}}_{W}^{j_{1}j_{2}} - \mathcal{H}_{W}^{j_{1}j_{2}}| \mathcal{F}_{n})^{2}\\
		=& E\left(\frac{1}{nr}\sum_{i=1}^{r}\frac{1}{\pi_{i}^{*}}w_{ij_{1}}^{*}w_{ij_{2}}^{*} - \frac{1}{n}\sum_{i=1}^{n}w_{ij_{1}}w_{ij_{2}}\bigg| \mathcal{F}_{n}\right)^{2}\\
		=& E\left(\frac{1}{nr}\sum_{i=1}^{r}\frac{1}{\pi_{i}^{*}}w_{ij_{1}}^{*}w_{ij_{2}}^{*}\bigg| \mathcal{F}_{n}\right)^{2} - \left(\frac{1}{n}\sum_{i=1}^{n}w_{ij_{1}}w_{ij_{2}}\right)^{2}\\
		=& \frac{1}{n^{2}r} \sum_{i=1}^{n}\pi_{i}\cdot\left(\frac{w_{ij_{1}}w_{ij_{2}}}{\pi_{i}}\right)^{2} - \frac{1}{n^{2}} \left(\sum_{i=1}^{n}w_{ij_{1}}w_{ij_{2}}\right)^{2}\\
		=& \frac{1}{n^{2}r} \sum_{i=1}^{n}\frac{(w_{ij_{1}}w_{ij_{2}})^{2}}{\pi_{i}} - \frac{1}{n^{2}} \left(\sum_{i=1}^{n}w_{ij_{1}}w_{ij_{2}}\right)^{2}\\
		\leq & \frac{1}{r} \max_{i=1,\ldots,n}(n\pi_{i})^{-1} \sum_{i=1}^{n}\frac{\|\mathbf{W}_{i}\|^{4}}{n}\\
		=&\frac{1}{r}O_{P}(1)O_{P}(1) = O_{P}(r^{-1}) .
	\end{align*}
	Using the Chebyshev inequality, for a sufficiently large $M$, we have
	\begin{align*}
		&P(\|\tilde{\mathcal{H}}_{W} - \mathcal{H}_{W} \| \geq r^{-1/2}M | \mathcal{F}_{n}) \\
		\leq& \frac{r E(\|\tilde{\mathcal{H}}_{W}\|^2 | \mathcal{F}_{n}) }{M^{2}}\\
		=&\frac{r \sum_{j_{1}=1}^{p}\sum_{j_{2}=1}^{p}E(\tilde{\mathcal{H}}_{W}^{j_{1}j_{2}} | \mathcal{F}_{n})^{2}}{M^{2}}\\
		=&\frac{O_{P}(1)}{M^{2}} \rightarrow 0, n,r \rightarrow \infty.
	\end{align*}
	Thus, equation (\ref{e3}) is proved.
	
	By (\ref{e3}) and Assumption \ref{js1}, we can obtain $\tilde{\mathcal{H}}_{W}^{-1} = O_{P|\mathcal{F}_{n}}(1)$. Therefore, combining (\ref{e1}), (\ref{e2}) and (\ref{e3}), we have
	$$
	\tilde{\boldsymbol{\beta}} - \hat{\boldsymbol{\beta}} = O_{P|\mathcal{F}_{n}}(r^{-1/2}).
	$$
	Then the theorem is proved.
\end{proof}

\bmhead{Proof of Theorem 2}
\begin{proof}
	Note that
	\begin{align*}
	\dot{\ell}^{*}(\hat{\boldsymbol{\beta}}) =& \frac{1}{r}\sum_{i=1}^{r} \left\{\frac{1}{n \pi_{i}^{*}} \left[-\mathbf{W}_{i}^{*} (y_{i}^{*} - \mathbf{W}_{i}^{*T}\hat{\boldsymbol{\beta}})\right] - \boldsymbol{\Sigma}_{uu}\hat{\boldsymbol{\beta}}\right\} \\
	 =& \frac{1}{r} \sum_{i=1}^{r} \boldsymbol{\xi}_{i},
    \end{align*}
	where $\boldsymbol{\xi}_{i} = \frac{1}{n \pi_{i}^{*}} \left[-\mathbf{W}_{i}^{*} (y_{i}^{*} - \mathbf{W}_{i}^{*T}\hat{\boldsymbol{\beta}})\right] - \boldsymbol{\Sigma}_{uu}\hat{\boldsymbol{\beta}}$ is an independent random vector.
	Then it can be directly obtained that
\begin{align}\label{e5}
			E(\boldsymbol{\xi}_{i}|\mathcal{F}_{n}) =& \frac{1}{n}\sum_{i=1}^{n} \left[-\mathbf{W}_{i} (y_{i} - \mathbf{W}_{i}^{T}\hat{\boldsymbol{\beta}})\right] - \boldsymbol{\Sigma}_{uu}\hat{\boldsymbol{\beta}} = \mathbf{0},\nonumber\\
			Var(\boldsymbol{\xi}_{i}|\mathcal{F}_{n}) =& \frac{1}{n^{2}}\sum_{i=1}^{n} \frac{\mathbf{W}_{i}\mathbf{W}_{i}^{T} (y_{i} - \mathbf{W}_{i}^{T}\hat{\boldsymbol{\beta}})^{2}}{\pi_{i}}  \nonumber\\
			&- (\boldsymbol{\Sigma}_{uu}\hat{\boldsymbol{\beta}})^{\otimes 2} = rV_{c}.
\end{align}
	According to Lemma {\ref{lem3}}, Assumptions \ref{js3} and \ref{js4}, we have
	\begin{align*}
		&\sum_{i=1}^{r} E\{\|r^{-1/2}\boldsymbol{\xi}_{i}\|^{2} I(\|r^{-1/2}\boldsymbol{\xi}_{i}\| > \varepsilon)|\mathcal{F}_{n}\}\\
		\leq & \frac{1}{\varepsilon^{\delta}r^{1+\delta/2}}\sum_{i=1}^{r} E\{\|\boldsymbol{\xi}_{i}\|^{2 + \delta}|\mathcal{F}_{n}\}\\
		\leq & \frac{2^{1+\delta}}{\varepsilon^{\delta}r^{1+\delta/2}} \sum_{i=1}^{r} \Bigg\{E\left[\left\|\frac{1}{n \pi_{i}^{*}} \mathbf{W}_{i}^{*} (y_{i}^{*} - \mathbf{W}_{i}^{*T}\hat{\boldsymbol{\beta}})\right\|^{2+\delta}\bigg|\mathcal{F}_{n}\right] \\
		&+ E[\|\boldsymbol{\Sigma}_{uu}\hat{\boldsymbol{\beta}}\|^{2+\delta}|\mathcal{F}_{n}]\Bigg\}\\
		= & \frac{2^{1+\delta}}{\varepsilon^{\delta}r^{\delta/2}} \left\{\sum_{i=1}^{n}\frac{ \|\mathbf{W}_{i}\|^{2+\delta} (y_{i} - \mathbf{W}_{i}^{T}\hat{\boldsymbol{\beta}})^{2+\delta}}{n^{2+\delta}\pi_{i}^{1+\delta}} + \|\boldsymbol{\Sigma}_{uu}\hat{\boldsymbol{\beta}}\|^{2+\delta}\right\}\\
		\leq & \frac{2^{1+\delta}}{\varepsilon^{\delta}r^{\delta/2}} \Bigg\{ \max_{i=1,\ldots,n}(n\pi_{i})^{-1-\delta} \sum_{i=1}^{n}\frac{(y_{i} - \mathbf{W}_{i}^{T}\hat{\boldsymbol{\beta}})^{2+\delta}\|\mathbf{W}_{i}\|^{2+\delta}}{n} \\
		&+ \|\boldsymbol{\Sigma}_{uu}\hat{\boldsymbol{\beta}}\|^{2+\delta}\Bigg\}\\
		=& O_{P}(r^{-\delta/2}).
	\end{align*}
	In the light of the Lindeberg-Feller central limit theorem \citep{Vaart1998}, it follows that
	\begin{equation}\label{e6}
		\left(\sum_{i=1}^{r} Var(\boldsymbol{\xi}_{i}|\mathcal{F}_{n})\right)^{-1/2}\sum_{i=1}^{r}\boldsymbol{\xi}_{i} = V_{c}^{-1/2} \dot{\ell}^{*}(\hat{\boldsymbol{\beta}}) \xrightarrow{d} N_{p}(\mathbf{0},I).
	\end{equation}
	By (\ref{e3}), we can obtain
	\begin{equation}\label{e7}
		\tilde{\mathcal{H}}_{W}^{-1} - \mathcal{H}_{W}^{-1} = -\mathcal{H}_{W}^{-1}(\tilde{\mathcal{H}}_{W} - \mathcal{H}_{W})\tilde{\mathcal{H}}_{W}^{-1} = O_{P|\mathcal{F}_{n}}(r^{-1/2}).
	\end{equation}
	By Assumption \ref{js1}, $\mathcal{H}_{W}$ converges to a positive definite matrix, then $\mathcal{H}_{W}^{-1} = O_{P}(1)$. And due to (\ref{e5}), we have
	\begin{equation}\label{e8}
		V= \mathcal{H}_{W}^{-1} V_{c} \mathcal{H}_{W}^{-1} = \frac{1}{r} \mathcal{H}_{W}^{-1} (rV_{c}) \mathcal{H}_{W}^{-1} = O_{P}(r^{-1}).
	\end{equation}
	Therefore, combining (\ref{e1}), (\ref{e7}) and (\ref{e8}), we have
	\begin{align*}
		&V^{-1/2}(\tilde{\boldsymbol{\beta}} - \hat{\boldsymbol{\beta}})\\
		=& -V^{-1/2} \tilde{\mathcal{H}}_{W}^{-1} \dot{\ell}^{*}(\hat{\boldsymbol{\beta}})\\
		=& -V^{-1/2} \mathcal{H}_{W}^{-1} \dot{\ell}^{*}(\hat{\boldsymbol{\beta}}) - V^{-1/2}(\tilde{\mathcal{H}}_{W}^{-1} - \mathcal{H}_{W}^{-1}) \dot{\ell}^{*}(\hat{\boldsymbol{\beta}})\\
		=& -V^{-1/2} \mathcal{H}_{W}^{-1}V_{c}^{1/2}V_{c}^{-1/2} \dot{\ell}^{*}(\hat{\boldsymbol{\beta}}) + O_{P|\mathcal{F}_{n}}(r^{-1/2}).
	\end{align*}
	Note that
	\begin{align*}&V^{-1/2} \mathcal{H}_{W}^{-1}V_{c}^{1/2}(V^{-1/2} \mathcal{H}_{W}^{-1}V_{c}^{1/2})^{T} \\
		=& V^{-1/2} \mathcal{H}_{W}^{-1}V_{c}^{1/2}V_{c}^{1/2}\mathcal{H}_{W}^{-1}V^{-1/2} = I, \end{align*}
	according to Slutsky's theorem \citep{Vaart1998}, we have
	\begin{equation}\label{e9}
		V^{-1/2}(\tilde{\boldsymbol{\beta}} - \hat{\boldsymbol{\beta}}) \xrightarrow{d} N_{p}(\mathbf{0},I), r , n \rightarrow \infty .
	\end{equation}
	Then the theorem is proved.
\end{proof}

\bmhead{Proof of Theorem 3}
\begin{proof}
	First of all, we prove (i). Let $V_{c}^{*} = \frac{1}{r n^{2}} \sum_{i=1}^{n} \frac{(y_{i} - \mathbf{W}_{i}^{T}\hat{\boldsymbol{\beta}})^{2}\mathbf{W}_{i}\mathbf{W}_{i}^{T}}{\pi_{i}}$.
	Note that $V_{c} = V_{c}^{*} -  \frac{1}{r}(\boldsymbol{\Sigma}_{uu}\boldsymbol{\beta})^{\otimes 2}$, the second item is not related to $\pi_{i}$. Therefore, to minimize $tr(V)$, we need to minimize $tr(\mathcal{H}_{W}^{-1} V_{c}^{*} \mathcal{H}_{W}^{-1})$.
	Then, we have
	\begin{align*}
		&tr(\mathcal{H}_{W}^{-1} V_{c}^{*} \mathcal{H}_{W}^{-1})\\
		=& \frac{1}{r n^{2}}\sum_{i=1}^{n} tr \left[\frac{(y_{i} - \mathbf{W}_{i}^{T}\hat{\boldsymbol{\beta}})^{2}\mathcal{H}_{W}^{-1} \mathbf{W}_{i} \mathbf{W}_{i}^{T} \mathcal{H}_{W}^{-1}}{\pi_{i}}\right]\\
		=& \frac{1}{r n^{2}}\sum_{i=1}^{n} \frac{(y_{i} - \mathbf{W}_{i}^{T}\hat{\boldsymbol{\beta}})^{2}\mathcal \|\mathcal{H}_{W}^{-1}\mathbf{W}_{i}\|^{2}}{\pi_{i}} \\
		=& \frac{1}{r n^{2}}\left(\sum_{j=1}^{n}\pi_{j}\right)\left[\sum_{i=1}^{n} \frac{(y_{i} - \mathbf{W}_{i}^{T}\hat{\boldsymbol{\beta}})^{2} \|\mathcal{H}_{W}^{-1}\mathbf{W}_{i}\|^{2}}{\pi_{i}}\right]\\
		\geq& \frac{1}{r n^{2}}\left[\sum_{i=1}^{n} |y_{i} - \mathbf{W}_{i}^{T}\hat{\boldsymbol{\beta}}|\cdot\|\mathcal{H}_{W}^{-1}\mathbf{W}_{i}\|\right]^{2}.
	\end{align*}
	The last inequality is the application of the Cauchy-Schwarz inequality. The equality holds if and only if $\pi_{i} \propto |y_{i} - \mathbf{W}_{i}^{T}\hat{\boldsymbol{\beta}}|\cdot\|\mathcal{H}_{W}^{-1}\mathbf{W}_{i}\|$. (ii) can be proved in the same way, so we omit the details here.
\end{proof}

\bmhead{Proof of Theorems 4 and 5}

Since the proofs of these two theorems are similar to those of Theorems \ref{theo1} and \ref{theo2}, we omit the proofs here.

\bmhead{Proof of Theorem 6}
\begin{proof}
	First, we consider the case where $m=1$. Note that
	\begin{equation}\label{eqs1}
		\begin{aligned}
			\check{\boldsymbol{\beta}} - \hat{\boldsymbol{\beta}}
			=& \left[\frac{1}{n} \sum_{i=1}^{n}\psi_{i}\mathbf{W}_{i}\mathbf{W}_{i}^{T} - \boldsymbol{\Sigma}_{uu}\right]^{-1} \\
			&\cdot \left[\frac{1}{n}\sum_{i=1}^{n} \psi_{i}\mathbf{W}_{i} (y_{i} - \mathbf{W}_{i}^{T}\hat{\boldsymbol{\beta}}) + \boldsymbol{\Sigma}_{uu}\hat{\boldsymbol{\beta}}\right]\\
			=& -(\check{\mathcal{H}}_{W})^{-1} \dot{L}^{*}(\hat{\boldsymbol{\beta}}),
		\end{aligned}
	\end{equation}
	where
	$$\check{\mathcal{H}}_{W} = \frac{1}{n} \sum_{i=1}^{n}\psi_{i}(\mathbf{W}_{i}\mathbf{W}_{i}^{T}) - \boldsymbol{\Sigma}_{uu},$$ $$\dot{L}^{*}(\boldsymbol{\beta}) = \frac{1}{n}\sum_{i=1}^{n} \psi_{i}\left[-\mathbf{W}_{i} (y_{i} - \mathbf{W}_{i}^{T}\boldsymbol{\beta})\right] - \boldsymbol{\Sigma}_{uu}\boldsymbol{\beta}.$$
	Therefore, it is only necessary to demonstrate that
	\begin{equation}\label{eqs2}
		\dot{L}^{*}(\hat{\boldsymbol{\beta}}) = O_{P|\mathcal{F}_{n}}(r^{-1/2}),
	\end{equation}
	and
	\begin{equation}\label{eqs3}
		\check{\mathcal{H}}_{W} - \mathcal{H}_{W} = O_{P|\mathcal{F}_{n}}(r^{-1/2}),
	\end{equation}
	where $\mathcal{H}_{W}= \frac{1}{n} \sum_{i=1}^{n}\mathbf{W}_{i}\mathbf{W}_{i}^{T} - \boldsymbol{\Sigma}_{uu}.$
	
	Note that
	$$E(\psi_{i}) = E(\mu_{i} \nu_{i}) = q_{n} \cdot \frac{1}{q_{n}} = 1,$$
	$$E(\psi_{i}^{2}) =\left[q_{n}(1-q_{n}) + q_{n}^{2}\right]\cdot\left(b_{n}^{2} + \frac{1}{q_{n}^{2}}\right) = q_{n} b_{n}^{2} + \frac{1}{q_{n}},$$
	and
	$$
	Var(\psi_{i}) = E(\psi_{i}^{2}) - [E(\psi_{i})]^{2} = q_{n} b_{n}^{2} + \frac{1}{q_{n}} - 1 = \frac{n a_{n}}{r} ,
	$$
	where $a_{n} = b_{n}^{2} q_{n}^{2} - q_{n} + 1.$
	According to Assumption \ref{js6}, as $q_{n}\rightarrow 0$, $\limsup\limits_{n\rightarrow \infty} a_{n} = \limsup\limits_{n\rightarrow \infty} q_{n}Var(\psi) = \limsup\limits_{n\rightarrow \infty} q_{n}(E(\psi_{i}^{2})-1) < \infty.$
	
	To prove (\ref{eqs2}), we can directly obtain
\begin{align*}
	&E(\dot{L}^{*}(\hat{\boldsymbol{\beta}})|\mathcal{F}_{n}) \\
		=& E\left\{\frac{1}{n}\sum_{i=1}^{n}\psi_{i}\left[-\mathbf{W}_{i} (y_{i} - \mathbf{W}_{i}^{T}\hat{\boldsymbol{\beta}})\right] - \boldsymbol{\Sigma}_{uu}\hat{\boldsymbol{\beta}} \bigg|\mathcal{F}_{n}\right\} \\
		=& \frac{1}{n}\sum_{i=1}^{n}\left[-\mathbf{W}_{i} (y_{i} - \mathbf{W}_{i}^{T}\hat{\boldsymbol{\beta}})\right] - \boldsymbol{\Sigma}_{uu}\hat{\boldsymbol{\beta}} = \mathbf{0}.
\end{align*}
	For the $j$-th element of $\dot{L}^{*}(\hat{\boldsymbol{\beta}})$, denote it as
	$$\dot{L}_{j}^{*}(\hat{\boldsymbol{\beta}}) = \frac{1}{n}\sum_{i=1}^{n}\psi_{i}[-w_{ij} (y_{i} - \mathbf{W}_{i}^{T}\hat{\boldsymbol{\beta}})] - (\boldsymbol{\Sigma}_{uu}\hat{\boldsymbol{\beta}})_{j}.
	$$
	By (\ref{ee1}) and Assumption \ref{js6}, we have
	\begin{align*}
		&Var(\dot{L}_{j}^{*}(\hat{\boldsymbol{\beta}})|\mathcal{F}_{n}) \\
		=& Var\left\{\frac{1}{n}\sum_{i=1}^{n}\psi_{i}\left[-w_{ij} (y_{i} - \mathbf{W}_{i}^{T}\hat{\boldsymbol{\beta}})\right] - (\boldsymbol{\Sigma}_{uu}\hat{\boldsymbol{\beta}})_{j}\bigg|\mathcal{F}_{n}\right\}\\
		=& \frac{1}{n^{2}}\frac{na_{n}}{r} \sum_{i=1}^{n}w_{ij}^{2} (y_{i} - \mathbf{W}_{i}^{T}\hat{\boldsymbol{\beta}})^{2}\\
		\leq& \frac{a_{n}}{r}\sum_{i=1}^{n} \frac{\|\mathbf{W}_{i}\|^{2}(y_{i} - \mathbf{W}_{i}^{T}\hat{\boldsymbol{\beta}})^{2}}{n}\\
		=& \frac{1}{r}O_{P}(1) = O_{P}(r^{-1}).
	\end{align*}
	From the Chebyshev inequality, for a sufficiently large $M$, we have
	\begin{align*}
		P(\|\dot{L}^{*}(\hat{\boldsymbol{\beta}})\| \geq r^{-1/2}M |\mathcal{F}_{n})
		&\leq \frac{r E(\|\dot{L}^{*}(\hat{\boldsymbol{\beta}})\|^{2}|\mathcal{F}_{n})}{M^{2}} \\
		&= \frac{r \sum_{j=1}^{p}E(\dot{L}_{j}^{*}(\hat{\boldsymbol{\beta}})|\mathcal{F}_{n})^{2}}{M^{2}} \\
		&= \frac{O_{P}(1)}{M^{2}}\rightarrow 0 , n, r \rightarrow \infty.
	\end{align*}
	Thus, equation (\ref{eqs2}) is proved.
	
	To prove (\ref{eqs3}), We can directly obtain
	$$E(\check{\mathcal{H}}_{W} | \mathcal{F}_{n}) = \mathcal{H}_{W} ,$$
	For any element $\check{\mathcal{H}}_{W}^{j_{1}j_{2}}$ of $\check{\mathcal{H}}_{W}$, $1\leq j_{1}, j_{2} \leq p$, by Assumptions \ref{js2} and \ref{js6}, we have
	\begin{align*}
		&Var(\check{\mathcal{H}}_{W}^{j_{1}j_{2}} | \mathcal{F}_{n}) \\
		=& Var\left[\frac{1}{n} \sum_{i=1}^{n}\psi_{i}(W_{ij_{1}}W_{ij_{2}}) - (\boldsymbol{\Sigma}_{uu})_{j_{1}j_{2}} \Big| \mathcal{F}_{n}\right]\\
		=& \frac{1}{n^{2}}\frac{na_{n}}{r} \sum_{i=1}^{n} (W_{ij_{1}}W_{ij_{2}})^{2} \\
		\leq& \frac{a_{n}}{r} \sum_{i=1}^{n}\frac{\|\mathbf{W}_{i}\|^{4}}{n} \\
		=& \frac{1}{r}O_{P}(1) = O_{P}(r^{-1}).
	\end{align*}
	From the Chebyshev inequality, for a sufficiently large $M$, we have
	\begin{align*}
		&P(\|\check{\mathcal{H}}_{W} - \mathcal{H}_{W} \| \geq r^{-1/2}M | \mathcal{F}_{n}) \\
		\leq& \frac{r E(\|\check{\mathcal{H}}_{W}\|^2 | \mathcal{F}_{n})}{M^{2}}\\
		=&\frac{r \sum_{j_{1}=1}^{p}\sum_{j_{2}=1}^{p}E(\check{\mathcal{H}}_{W}^{j_{1}j_{2}} | \mathcal{F}_{n})^{2}}{M^{2}}\\
		=&\frac{O_{P}(1)}{M^{2}} \rightarrow 0, n,r \rightarrow \infty.
	\end{align*}
	Thus, equation (\ref{eqs3}) is proved.
	
	By (\ref{eqs3}) and Assumption \ref{js1}, we have $\check{\mathcal{H}}_{W}^{-1} = O_{P|\mathcal{F}_{n}}(1)$. Therefore, combining (\ref{eqs1}), (\ref{eqs2}) and (\ref{eqs3}), we have
	\begin{equation}\label{eqs4}
		\check{\boldsymbol{\beta}} - \hat{\boldsymbol{\beta}} = O_{P|\mathcal{F}_{n}}(r^{-1/2}).
	\end{equation}
	
	As $m >1$ , we have $\check{\boldsymbol{\beta}}^{(m)} = \frac{1}{m} \sum_{k=1}^{m}\check{\boldsymbol{\beta}}_{k}.$
	Then according to the weak law of large numbers, it follows that
	\begin{align*}
		\check{\boldsymbol{\beta}}^{(m)} - \hat{\boldsymbol{\beta}} = \frac{1}{m} \sum_{k=1}^{m}\check{\boldsymbol{\beta}}_{k} - \hat{\boldsymbol{\beta}} =& \frac{1}{m} \sum_{k=1}^{m}(\check{\boldsymbol{\beta}}_{k} - \hat{\boldsymbol{\beta}}) \\
		=& O_{P|\mathcal{F}_{n}}((mr)^{-1/2}).
		\end{align*}
	Then the theorem is proved.
\end{proof}

\bmhead{Proof of Theorem 7}
\begin{proof}
	Firstly, we prove the case where $m=1$.
	Because
\begin{align*}
	\dot{L}^{*}(\hat{\boldsymbol{\beta}}) &= \frac{1}{n}\sum_{i=1}^{n}\left\{\psi_{i}\left[-\mathbf{W}_{i} (y_{i} - \mathbf{W}_{i}^{T}\hat{\boldsymbol{\beta}})\right] - \boldsymbol{\Sigma}_{uu}\hat{\boldsymbol{\beta}}\right\}  \\
	&= \frac{1}{\sqrt{r}} \sum_{i=1}^{n} \boldsymbol{\eta}_{i},
\end{align*}
	where $\boldsymbol{\eta}_{i} = \frac{\sqrt{r}}{n}\left\{\psi_{i}\left[-\mathbf{W}_{i} (y_{i} - \mathbf{W}_{i}^{T}\hat{\boldsymbol{\beta}})\right] - \boldsymbol{\Sigma}_{uu}\hat{\boldsymbol{\beta}}\right\}$ is an independent random vector.
	Note that
	\begin{align*}
		E(\boldsymbol{\eta}_{i}|\mathcal{F}_{n}) &= \frac{\sqrt{r}}{n}\left[-\mathbf{W}_{i} (y_{i} - \mathbf{W}_{i}^{T}\hat{\boldsymbol{\beta}}) - \boldsymbol{\Sigma}_{uu}\hat{\boldsymbol{\beta}}\right],\\
		Var(\boldsymbol{\eta}_{i}|\mathcal{F}_{n}) &= \frac{r}{n^{2}}\frac{n a_{n}}{r}\left[-\mathbf{W}_{i} (y_{i} - \mathbf{W}_{i}^{T}\hat{\boldsymbol{\beta}})\right]^{\otimes 2} \\
		&= \frac{a_{n}}{n}\mathbf{W}_{i}\mathbf{W}_{i}^{T} (y_{i} - \mathbf{W}_{i}^{T}\hat{\boldsymbol{\beta}})^{2}.
	\end{align*}
	Then by using Assumptions \ref{js2} and \ref{js6}, we obtain
\begin{align}
	\sum_{i=1}^{n} E(\boldsymbol{\eta}_{i}|\mathcal{F}_{n}) =& \frac{\sqrt{r}}{n}\sum_{i=1}^{n}\left[-\mathbf{W}_{i} (y_{i} - \mathbf{W}_{i}^{T}\hat{\boldsymbol{\beta}}) - \boldsymbol{\Sigma}_{uu}\hat{\boldsymbol{\beta}}\right] \nonumber\\
	=& \mathbf{0},\nonumber\\
	\sum_{i=1}^{n} Var(\boldsymbol{\eta}_{i}|\mathcal{F}_{n}) &= \frac{a_{n}}{n}\sum_{i=1}^{n} \mathbf{W}_{i}\mathbf{W}_{i}^{T} (y_{i} - \mathbf{W}_{i}^{T}\hat{\boldsymbol{\beta}})^{2} \nonumber\\
		&=a_{n} \Sigma_{c}.\label{eqs5}
\end{align}
	According to Lemma {\ref{lem3}}, Assumptions \ref{js4} and \ref{js6}, we have
	\begin{align*}
		&\sum_{i=1}^{n} E\{\|\boldsymbol{\eta}_{i}\|^{2} I(\|\boldsymbol{\eta}_{i}\| > \varepsilon)|\mathcal{F}_{n}\}\\
		\leq & \varepsilon^{-\alpha}\sum_{i=1}^{n} E\{\|\boldsymbol{\eta}_{i}\|^{2 + \alpha}|\mathcal{F}_{n}\}\\
		= & \frac{r^{1+ \frac{\alpha}{2}}}{\varepsilon^{\alpha} n^{2+\alpha}}\sum_{i=1}^{n} E\left\{\left\|\psi_{i}\mathbf{W}_{i} (y_{i} - \mathbf{W}_{i}^{T}\hat{\boldsymbol{\beta}}) + \boldsymbol{\Sigma}_{uu}\hat{\boldsymbol{\beta}}\right\|^{2+\alpha}\bigg|\mathcal{F}_{n}\right\}\\
		\leq & \frac{r^{1+\frac{\alpha}{2}} 2^{1+\alpha}}{\varepsilon^{\alpha} n^{2+\alpha}} \sum_{i=1}^{n} \Bigg\{E\left[\left\|\psi_{i}\mathbf{W}_{i} (y_{i} - \mathbf{W}_{i}^{T}\hat{\boldsymbol{\beta}})\right\|^{2+\alpha}\bigg|\mathcal{F}_{n}\right]\\ &+E\left[\left\|\boldsymbol{\Sigma}_{uu}\hat{\boldsymbol{\beta}}\right\|^{2+\alpha}\bigg|\mathcal{F}_{n}\right]\Bigg\}\\
		=& \frac{r^{1+\frac{\alpha}{2}} 2^{1+\alpha}}{\varepsilon^{\alpha} n^{1+\alpha}} \Bigg\{E(\psi)^{2+\alpha}\frac{1}{n} \sum_{i=1}^{n} \|\mathbf{W}_{i}\|^{2+\alpha} (y_{i} - \mathbf{W}_{i}^{T}\hat{\boldsymbol{\beta}})^{2+\alpha} \\&+\left\|\boldsymbol{\Sigma}_{uu}\hat{\boldsymbol{\beta}}\right\|^{2+\alpha} \Bigg\} \\
		= & O_{P}(r^{-\alpha/2}).
	\end{align*}
	Therefore, the Lindeberg-Feller condition is satisfied. According to the Lindeberg-Feller central limit theorem, we have
	\begin{align}\label{eqs6}
		&\left(\sum_{i=1}^{n} Var(\boldsymbol{\eta}_{i}|\mathcal{F}_{n})\right)^{-1/2}\sum_{i=1}^{n}\boldsymbol{\eta}_{i} \nonumber\\
		=& \sqrt{\frac{r}{a_{n}}} \Sigma_{c}^{-1/2} \dot{L}^{*}(\hat{\boldsymbol{\beta}}) \xrightarrow{d} N_{p}(\mathbf{0},I).
	\end{align}
	By (\ref{eqs3}), we have
	\begin{equation}\label{eqs7}
		\check{\mathcal{H}}_{W}^{-1} - \mathcal{H}_{W}^{-1} = -\mathcal{H}_{W}^{-1}(\check{\mathcal{H}}_{W} - \mathcal{H}_{W})\check{\mathcal{H}}_{W}^{-1} = O_{P|\mathcal{F}_{n}}(r^{-1/2}).
	\end{equation}
	By Assumption \ref{js1}, $\mathcal{H}_{W}$ converges to a positive definite matrix, then $\mathcal{H}_{W}^{-1} = O_{P}(1)$.
	And due to (\ref{eqs5}), we obtain
	\begin{equation}\label{eqs8}
		\Sigma= \mathcal{H}_{W}^{-1} \Sigma_{c} \mathcal{H}_{W}^{-1} = O_{P}(1).
	\end{equation}
	Therefore, combining(\ref{eqs1}), (\ref{eqs7}) and (\ref{eqs8}), we have
	\begin{align*}
		&\sqrt{\frac{r}{a_{n}}}\Sigma^{-1/2}(\check{\boldsymbol{\beta}} - \hat{\boldsymbol{\beta}})\\
		 =& - \sqrt{\frac{r}{a_{n}}}\Sigma^{-1/2} \check{\mathcal{H}}_{W}^{-1} \dot{L}^{*}(\hat{\boldsymbol{\beta}})\\
		=& -\sqrt{\frac{r}{a_{n}}}\Sigma^{-1/2} \mathcal{H}_{W}^{-1} \dot{L}^{*}(\hat{\boldsymbol{\beta}}) \\
		&- \sqrt{\frac{r}{a_{n}}}\Sigma^{-1/2}(\check{\mathcal{H}}_{W}^{-1} -\mathcal{H}_{W}^{-1}) \dot{L}^{*}(\hat{\boldsymbol{\beta}})\\
		=& -\sqrt{\frac{r}{a_{n}}}{\Sigma}^{-1/2} \mathcal{H}_{W}^{-1}\Sigma_{c}^{1/2}\Sigma_{c}^{-1/2} \dot{L}^{*}(\hat{\boldsymbol{\beta}}) + O_{P|\mathcal{F}_{n}}(r^{-1/2}).
	\end{align*}
	Note that
		\begin{align*}
			&\Sigma^{-1/2} \mathcal{H}_{W}^{-1}\Sigma_{c}^{1/2}(\Sigma^{-1/2} \mathcal{H}_{W}^{-1}\Sigma_{c}^{1/2})^{T} \\
			=& \Sigma^{-1/2} \mathcal{H}_{W}^{-1}\Sigma_{c}^{1/2}\Sigma_{c}^{1/2}\mathcal{H}_{W}^{-1}\Sigma^{-1/2} = I,
		\end{align*}
	From the Slutsky theorem and (\ref{eqs6}), we have
	\begin{equation*}
		\Sigma^{-1/2} \sqrt{r/a_{n}}(\check{\boldsymbol{\beta}} - \hat{\boldsymbol{\beta}}) \xrightarrow{d} N_{p}(\mathbf{0},I), r , n \rightarrow \infty .
	\end{equation*}
	
	As $m > 1$, we have $\check{\boldsymbol{\beta}}^{(m)} = \frac{1}{m} \sum_{k=1}^{m}\check{\boldsymbol{\beta}}_{k}.$
	Because $\psi_{k,i}$ are independent of each other when $\mathcal{F}_{n}$ is known, then $\check{\boldsymbol{\beta}}_{k}$ are independent of each other for $k = 1, \cdots, m$.
	From the multivariate sampling distribution theorem, we can obtain that
	\begin{align*}
		&\sqrt{\frac{rm}{a_{n}}}(\check{\boldsymbol{\beta}}^{(m)} - \hat{\boldsymbol{\beta}})\\
		=& \frac{1}{\sqrt{m}} \sum_{k=1}^{m}\sqrt{\frac{r}{a_{n}}}(\check{\boldsymbol{\beta}}_{k} - \hat{\boldsymbol{\beta}})
		\xrightarrow{d} N_{p}(\mathbf{0},\Sigma), r , n \rightarrow \infty.
	\end{align*}
	Then the theorem is proved.
\end{proof}

\end{appendices}


\bibliography{sn-bibliography}

\end{document}